\documentclass[11pt]{ps}

\usepackage{graphics,graphicx}

\def\a{\alpha }

\def\sJ{\mathcal{J}}

\numberwithin{equation}{section}

\begin{document}

\title{On a Szeg\"o type limit theorem,\\ the H\"older-Young-Brascamp-Lieb
inequality, \\and the asymptotic theory of
integrals and quadratic forms of stationary fields}
\author{Florin Avram}\address{Dept. de Math., Universite de Pau, France, E-mail: Florin.Avram@univ-Pau.fr}
\author{Nikolai Leonenko} \address{
Cardiff School of Mathematics, Cardiff University, Senghennydd
Road, Cardiff, CF24 4AG, UK, \\ Email: LeonenkoN@Cardiff.ac.uk}
\author{Ludmila Sakhno}\address{Dept. of Probability Theory and Mathematical Statistics, Kyiv National Taras
Shevchenko University, Ukraine, and Cardiff School of Mathematics, Cardiff University, Senghennydd
Road, Cardiff, CF24 4AG, UK, Email:lms@univ.kiev.ua}
\thanks{Partly supported by the grant of London Mathematical Society 2632, the EPSRC
grant RCMT 119 and by the Welsh Institute of \ Mathematics and Computational
Sciences.}

\begin{abstract}
\noindent Many statistical applications require establishing
central limit theorems for sums/integrals $S_T(h)=\int_{t \in I_T}
h (X_t) dt$ or for quadratic forms $Q_T(h)=\int_{t,s \in I_T}
\hat{b}(t-s) h (X_t, X_s) ds dt$, where $X_t$ is a stationary
process. A particularly important case is that of Appell
polynomials $h (X_t)=P_{m} (X_t)$, $h (X_t,X_s) = P_{m,n} (X_t,
X_s)$, since the ``Appell expansion rank" determines typically the
type of central limit theorem satisfied by the functionals
$S_T(h), Q_T(h)$.
We review and extend here to multidimensional 
indices, along lines conjectured in \cite{AT06}, a functional
analysis approach to this problem proposed by Avram and Brown
(1989), based on the method of cumulants and on integrability
assumptions in the spectral domain; several applications are
presented as well.
\end{abstract}

\begin{resume}
Nous consid\'erons ici des  th\'eor\'emes de limite central  pour
des sommes/integrales $S_T(h)=\int_{t \in I_T} h (X_t) dt,$ et pour
des formes quadratiques $Q_T(h)=\int_{t,s \in I_T} \hat{b}(t-s) h
(X_t, X_s) ds dt$, o\`u $X_t$ est un processus stationnaire. Un cas
particuli\`erement important est celui des polyn\^omes d'Appell
$h(X_t)=P_{m} (X_t)$, $h (X_t,X_s) = P_{m,n} (X_t, X_s)$.
Pour ce probl\`eme, nous     g\'en\'eralisons  ici au cas des
indices multidimensionnels une approche
propos\'ee par
Avram et Brown (1989), bas\'ee sur la m\'ethode des cumulants et sur
des hypoth\`eses d'integrabilit\'e dans le  domaine spectral.
Plusieurs applications illustrent la versatilit\'e de l'approche.
\end{resume}

\maketitle

 \tableofcontents

\section*{Introduction}

\label{s:intro}

\textbf{Model.} For a unified treatment of the discrete and the continuous
(multi-dimensional) setups, we assume that $X_{t},t\in I$ is a real
stationary random field, where $I$ denotes a measurable group endowed with
its Haar measure, and $\int_{t\in I}\nu (dt)$ is integral with respect to
the Haar measure. Depending on the continuous/discrete setup, $I$ will be
either $\mathbb{R}^{d}$ with Lebesgue measure, or $\mathbb{Z}^{d}$ with the
counting measure. Discrete sums might be written either as integrals (in the
statement of theorems), or in traditional sum notation (in the expository
part). For continuous case we assume that $X_{t},t\in \mathbb{R}^{d}$ is a
measurable random field.

We will assume throughout the existence of all order cumulants $%
c_{k}(t_{1},t_{2},...,t_{k})$ for our stationary in the strict sense random
field $X_{t}$, which are supposed to be representable as Fourier transforms
of \textquotedblleft cumulant spectral densities" $f_{k}(\lambda
_{1},...,\lambda _{k-1})\in L_{1}$, i.e:
\begin{equation*}
c_{k}(t_{1},t_{2},...,t_{k})=c_{k}(t_{1}-t_{k},..,t_{k-1}-t_{k},0)=
\end{equation*}%
\begin{equation*}
=\int_{\lambda _{1},...,\lambda _{k-1}\in S}e^{i\sum_{j=1}^{k-1}\lambda
_{j}(t_{j}-t_{k})}f_{k}(\lambda _{1},...,\lambda _{k-1})\,\mu (d\lambda
_{1})...\mu (d\lambda _{k-1}).
\end{equation*}%
Throughout, $S$ will denote the \textquotedblleft spectral" space of
discrete/continuous processes and $\mu (d\lambda )$ will denote the
associated Haar measure, i.e. Lebesgue measure normalized to unity on $[-\pi
,\pi ]^{d}$ and Lebesgue measure on $\mathbb{R}^{d}$, respectively. The
functions $f_{k}(\lambda _{1},...,\lambda _{k-1})$ are symmetric and
generally complex valued for real field $X_{t}$.

The random field $X_t $ is observed on a sequence $I_T$ of increasing
dilations of a finite convex domain $I_1$, i.e.
\begin{equation*}
I_T= T \; I_1, \quad T \to \infty .
\end{equation*}

Correspondingly, we will consider linear and bilinear forms $S_{T}$ and $%
Q_{T}$, where summation/integration is performed over domains $I_{T}$.

In the discrete-time case, the cases $I_{T}=[1,T]^{d}$, $T\in \mathbb{Z}_{+}$
(in keeping with tradition) or $I_{T}=[-T/2,T/2]^{d}$, $T\in 2\mathbb{Z}_{+}$
will be assumed. In the continuous case, we focus on rectangles $%
I_{T}=\{t\in \mathbb{R}^{d}:-T/2\leq t_{i}\leq T/2,i=1,...,d\}$.

Later on we will see that the choice of a specific set $I_{1}$ leads (when
working in the spectral domain) to using an associated Dirichlet type kernel
\begin{equation}
\Delta _{T}(\lambda )=\int_{t\in I_{T}}e^{it\lambda }\nu (dt)  \label{ker}
\end{equation}%
and a multivariate Fej\'{e}r kernel (\ref{Fk}). Explicit well-known formulas
are available when $I_{1}$ is a rectangle or a ball, both for the discrete
and continuous case. We will work in the paper with rectangular domains,
however, some extensions are possible for balls by replacing corresponding
estimates (see Appendix B).

\textbf{Motivation.} Via the work of Szego, Schur, Wiener and Kolmogorov,
the study of stationary one-dimensional discrete time series, defined by
choosing $I=\mathbb{Z}$, has been well developed and became tightly
interwoven with several branches of mathematics, bringing forth jewels like
for example the Wiener-Kolmogorov formula identifying the variance of the
prediction error with respect to the past with the integral of the log of
the spectral density. The extension to the continuous time case $I=\mathbb{R}
$, provided by Krein, exemplifies the subtle challenges brought by modifying
the nature of the index set.

The convenience of time series comes largely from the FARIMA family of
parametric models, defined as solutions of equations
\begin{equation*}
\phi (B)(1-B)^{d}X_{t}=\theta (B)\xi _{t}
\end{equation*}%
where $B$ is the operator of backward translation in time, $\phi (B),\theta
(B)$ are polynomials, $d$ is a real number and $\xi _{t}$ is white noise
\cite{hurst:1951}, \cite{Granger}, \cite{Hosking}, \cite{beran:1994B}, \cite%
{willinger:taqqu:teverovsky:1999}. Using the FARIMA family of models, one
attempts, via an extension of the Box-Jenkins methodology, to estimate the
parameter $d$ and the coefficients of the polynomials $\phi (B),\theta (B)$
such that the residuals $\xi _{t}$ display white noise behavior (and hence
may be safely discarded for prediction purposes).

An extension of this approach to continuous time and to multi-parameter
processes (spatial statistics) has been long time missing. Only recently,
the FICARMA (\cite{brockwell}, \cite{AHL}) and the Riesz-Bessel families of
processes (which essentially replace the operator $B$ by the Laplacian - see
Appendix A), have allowed a similar approach for processes with continuous
and multidimensional indices (at least in the rotationally invariant case).

These examples illustrate the point that extension of ARIMA-type models to
continuous and multidimensional indices is an important challenge.

\textbf{Central limit theorems for quadratic forms.} Asymptotic statistical
theory, and in particular, estimation of the parameters of FARIMA and
Riesz-Bessel processes, requires often establishing central limit theorems
concerning
\begin{eqnarray*}
\text{sums} &&\text{ \ }S_{T}=S_{T}(h)=\sum_{i=1}^{T}h(X_{i}) \\
\text{and bilinear forms \ } &&Q_{T}=Q_{T}(h)=\sum_{i=1}^{T}\sum_{j=1}^{T}%
\hat{b}(i-j)\;h(X_{i},X_{j}).
\end{eqnarray*}%
of functions of stationary sequences $X_{i}$ (and their generalizations to
continuous and multidimensional indicies).

\textbf{Notes:} 1) The kernel of the quadratic form will be denoted by $\hat{%
b}(t)$, so that we may reserve $b(\lambda)$ for its Fourier transform.

2) The random field $X_{t},$ $t\in I_{T},$ will be allowed to have
short-range or long-range dependence (that is, summable or non-summable
correlations).

A particularly important case is that of Appell polynomials\footnote{%
For the definition of Appell polynomials see, for instance, Avram and Taqqu
(1987), or Giraitis and Surgailis (1986).} $h(X_{t})=P_{m}(X_{t}),$ $%
h(X_{t},X_{s})=P_{m,n}(X_{t},X_{s})$ associated to the distribution of $%
X_{t} $, which are the building block of the so called \textquotedblleft
chaos/Fock expansions". Two main cases were distinguished, depending on
whether the limit is Gaussian or not (the latter case being referred to as a
non-central limit theorem).

This line of research, initiated by Dobrushin and Major \cite%
{dobrushin:major:1979} and Taqqu \cite{taqqu:1979} in the Gaussian case (see
also Ivanov and Leonenko (1989) for Gaussian continuous case), by Giraitis
and Surgailis \cite{Gir}, \cite{giraitis:surgailis:1986}, and by Avram and
Taqqu \cite{avram:taqqu:1987} in the linear case, continues to be of
interest today \cite{Gin}, \cite{giraitis:taqqu:1997}, \cite{gtw}, \cite{NP},
\cite{PT},\cite{Sol} to name only a few papers.

Our interest here is in obtaining extensions to continuous and
multi-parameter processes of the central limit theorem for sums and
quadratic form, obtained in the case of discrete time series by Breuer and
Major \cite{BM}, by the method of moments.

\textbf{Some analytic tools.} A key unifying role in our story is played by
the multivariate Fej\'{e}r kernels:
\begin{equation}
\Phi _{T}^{\ast }(u_{1},...,u_{n-1})=\frac{1}{\left( 2\pi \mu (I_{1})\right)
^{\left( n-1\right) d}T^{d}}\Delta
_{T}(-\sum_{1}^{n-1}u_{e})\prod_{e=1}^{n-1}\Delta _{T}(u_{e}),  \label{Fk}
\end{equation}%
and their kernel property: the fact that \textbf{when $T\rightarrow \infty $%
, the multivariate Fej\'{e}r kernel }$\Phi _{T}^{\ast }$ \textbf{convergence
weakly to a $\delta $ measure}:

\begin{lmm}
\label{lwc} \textbf{The kernel property:} For any continuous bounded
function $C(u_{1},..,u_{n-1})$, it holds that:
\begin{equation*}
\lim_{T\rightarrow \infty }\int_{S^{n-1}}C(u_{1},..,u_{n-1})\Phi _{T}^{\ast
}(u_{1},...,u_{n-1})\prod_{i=1}^{n-1} \mu(du_{i})=C(0,...,0)
\end{equation*}
\end{lmm}

\textbf{Proof:} For the discrete onedimensional case we refer for example to
\cite{AB}, and for the continuous multidimensional case, with $I_{1}$ a
rectangle, to Bentkus \cite{Ben}, \cite{Ben1} or \cite{ALS1}, Proposition 1.%
\footnote{%
This convergence of measures may also be derived as a consequence of the H%
\"{o}lder-Young-Brascamp-Lieb \ inequality (see Theorem \ref{t:HYB}), using
estimates of the form
\begin{equation*}
||{\Delta _{T}}||_{s_{v}^{-1}}\leq k(s_{v})T^{d(1-s_{v})}
\end{equation*}%
with optimally chosen $s_{v}$, $v=1,\cdots ,V$.}.

Developing some limit theory for multivariate Fej\'{e}r kernels was the key
point in several papers \cite{AB}, \cite{a}, \cite{AT}, \cite{af} which
generalized the Breuer and Major central limit theorem \cite{BM}. The papers
above introduced a new mathematical object to be called \textquotedblleft
\textbf{Fej\'{e}r graph kernels}" -- see (\ref{matF2}) in section \ref{s:kgi}%
, which captures the common structure of several cumulant computations.
Replacing the cyclic graph encountered in the specific case of quadratic
forms in Gaussian random variables by an arbitrary graph, these papers
reduce the central limit theorem for a large class of problems involving
Appell polynomials in Gaussian or moving average summands to an application
of three analytical tools:

\begin{enumerate}
\item Identifying the graphs involved by applying the well-known \textbf{%
diagram formula} for computing moments/cumulants of Wick products -- see
section \ref{s:cum}.

\item Applying a generalization of a \textbf{Grenander-Szeg\"{o} theorem on
the trace of products of Toeplitz matrices} to the {Fej\'{e}r graph integrals%
} -- see section \ref{s:kgi}, to obtain the asymptotic variance. This
theorem is valid under some general integrability assumptions furnished by
the H\"{o}lder-Young-Brascamp-Lieb inequality.

\item The resolution of certain combinatorial graph optimization problems,
specific to each application, which clarify the geometric structure of the
polytope of valid integrability exponents on the functions involved
(spectral density, kernel of the quadratic form, etc).
\end{enumerate}

Here, we observe that a similar approach works in the multidimensional and
continuous indices case. More precisely, the only changes are a) the
normalizations, which change from $T$ to $T^{d}$, and b) the condition for
the validity of the H\"{o}lder-Young-Brascamp-Lieb inequality (see Appendix
C) in the continuous case. Therefore, the previously obtained central limit
theorems continue to hold in the multidimensional case, including continuous
indices, after simply adjusting the normalizations and integrability
conditions.

\textbf{Statistical applications.} The two cases most easy to study are that
of Gaussian and linear processes. In the applications Section \ref{s:app} we
will work assuming that $X_{t}$ is a linear process (see (\ref{nn1}) or (\ref%
{nn3}) below). 
This assumption has the advantage of implying a product representation of
the cumulant spectral densities -- see for example Theorem 2.1 of \cite{AHL}%
. Namely, for the cumulants we get the explicit formula
\begin{equation*}
c_{k}(t_{1},...,t_{k})=d_{k}\int_{s\in I}\prod_{j=1}^{k}\hat{a}%
(t_{j}-s)\;\nu (ds),
\end{equation*}%
and in the spectral domain, we get
\begin{equation}
f_{k}(\lambda _{1},...,\lambda _{k-1})=d_{k}\;a(-\sum_{i=1}^{k-1}\lambda
_{i})\;\prod_{i=1}^{k-1}a(\lambda _{i})=d_{k}\;\prod_{i=1}^{k}a(\lambda
_{i})\delta (\sum_{j=1}^{k}\lambda _{j})  \label{lindens}
\end{equation}
(the meaning of parameters $d_{k}$ and a function $a(\lambda )$ is clarified
in section \ref{ss:lin}).

For $k=2$, we will denote the spectral density by $f(\lambda)=
f_2(\lambda)=d_2 a(\lambda) a(-\lambda)$.

\begin{dfntn}
Let
\begin{eqnarray}
{{\mathbf{L}}}_{p} (d \mu )=
\begin{cases}
{L}_{p} (d \mu ) & \text{ if } p< \infty , \\
C & \text{ if } p= \infty .%
\end{cases}%
\end{eqnarray}
denote the closure of the functions in $L_p$ which are \textbf{continuous,
bounded and of bounded support}, under the $L_p$ norm.
\end{dfntn}

\textbf{Note: } In the torus case, this space intervenes in the proof of
theorem \ref{Sz}, which is established first for complex exponentials, and
extended then to the Banach space of functions which may be approximated
arbitrarily close in $L_p(d \mu )$ sense by linear combinations of complex
exponentials, endowed with the $L_p$ norm.


Considering bilinear forms $Q_T$ we will work under integrability assumption:

\textbf{Assumption A:}
\begin{equation*}
f \in {\mathbf{L}}_{p_1}(d \mu), b \in {\mathbf{L}}_{p_2}(d \mu), \quad
1\leq p_i \leq \infty, i=1,2
\end{equation*}

\textbf{Note:} 1) While a general stationary model is parameterized by a
sequence of functions $f_{k}(\lambda _{1},...,\lambda _{k-1})$, $k\in
\mathbb{N}$, the linear model (\ref{lindens}) is considerably simpler, being
parametrized by a single function $a(\lambda )$.

2) We expect all our results may be formulated directly in terms of
characteristics of the field $X_{t}$, which suggests that the moving average
assumption is probably unnecessary; indeed, more general results which make
direct assumptions that functions $f_{k}(\lambda _{1},...,\lambda _{k-1})$
belong to some special $L_{p}$-type spaces, have been obtained in certain
cases -- see, for example, \cite{af}.

\textbf{Contents.} We present a warmup example involving quadratic Gaussian
forms in Section \ref{s:ex}. The results here are closely connected to those
of the paper by Ginovian and Sahakyan (2007). We define the concept of
Fej\'er graph and matroid integrals in Section \ref{s:kgi}. We will consider
here only the first case (i.e. graphic matroids associated to the incidence
matrix of a graph).

Some limit theory (of Grenander-Szeg\"{o} type) for Fej\'{e}r graph
integrals is reviewed and extended to the continuous case in Section \ref%
{s:GS}. Various estimates concerning kernels are collected in Appendix \ref%
{s:ker}, and a particular case of the H\"{o}lder-Young-Brascamp-Lieb
inequality required here is presented in Appendix \ref{s:HYBu}.

In Section \ref{s:app} we introduce the linear model (which extends the
Gaussian model) and develop several applications. Note here the existence of
a different approach, due to Peligrad and Utev (2006), who established the
central limit theorem for linear processes with discrete time and dependent
innovations including martingale and mixingale type assumptions (see also
the references therein for this line of investigation).

To make the paper self-contained we supply in Appendices the material we
refer to in the main part of the paper.

\section{An example: the central limit theorem for Gaussian bilinear forms
\label{s:ex}}

We present first our method in the simplest case of symmetric bilinear forms
$Q_{T}=Q_{T}^{(1,1)}=Q_{T}(P_{1,1})$ in stationary Gaussian fields $X_{t}$,
with covariances $r({t-s})$, $t,s\in I$ and spectral density $f(\lambda )$
(note that $S_{T}^{(1)}=\int_{t\in I_{T}}X_{t}\nu (dt)$ is \textquotedblleft
too simple" for our purpose, since it is already Gaussian and and its $k$-th
order cumulants $\chi _{k}(S_{T}^{(1)})=0,\;\forall k\neq 2$). The
presentation follows \cite{a} for the discrete case and \cite{Gin} for the
continuous case, except that we clarify the point that the previous results
hold in any dimension $d$.

To obtain the central limit theorem for $T^{-d/2}Q_{T}^{(1,1)}$ by the
method of cumulants it is enough to show that:
\begin{equation*}
\lim_{T\rightarrow \infty }\chi _{2}\left( \frac{Q_{T}^{(1,1)}}{T^{d/2}}%
\right) \text{ is finite, and }\quad \lim_{T\rightarrow \infty }\chi
_{k}\left( \frac{Q_{T}^{(1,1)}}{T^{d/2}}\right) =0,\;\forall k\geq 3.
\end{equation*}

A direct computation based on multilinearity 
yields the cumulants of $Q_T^{(1,1)}$:

\begin{equation}
\chi _{k}=\chi (Q_{T},...,Q_{T})=2^{k-1}(k-1)!\ \mathrm{Tr}%
[(T_{T}(b)T_{T}(f))^{k}].  \label{cumQ}
\end{equation}%
\noindent Here, $\mathrm{Tr}$ denotes the trace and
\begin{equation*}
T_{T}(b)=(\hat{b}(t{-s}),\ t,s\in I_{T}),\quad \quad \quad T_{T}(f)=({r}%
(t-s),\ t,s\in I_{T})
\end{equation*}%
denote Toeplitz matrices (with multidimensional indices) of dimension $%
T^{d}\times T^{d}$ in the discrete case and truncated Toeplitz-type
operators in the continuous case \footnote{%
Recall that in continuous case, the truncated Toeplitz-type operator
generated by a function $\hat{f}\in L_{\infty },$ is defined for $u\in L_{2%
\text{ }}$ as follows
\begin{equation*}
T_{T}({f})u(t)=\int\limits_{I_{T}}\hat{f}(t-s)u(s)\nu (ds).
\end{equation*}%
}.

While the cumulants in (\ref{cumQ}) may be expressed using powers of two
Toeplitz matrices, it turns out more convenient in fact to consider more
general products with all terms potentially different (taking advantage thus
of multi-linearity).

Suppose therefore given a set ${f}_e(\lambda): S \to \mathbb{R},\
e=1,\ldots,n$ of ``symbols" associated to the set of Toeplitz operators,
where $(S, d \mu)$ denotes either $\mathbb{R}^d$ with Lebesgue measure, or
the torus $[-\pi,\pi]^d$ with normalized Lebesgue measure. Assume the
symbols satisfy integrability conditions
\begin{equation}  \label{e:Lp1}
f_e \in {\mathbf{L}}_{p_e} (S, d \mu ),\ 1\leq p_e\leq \infty.
\end{equation}

Let $\hat{f}_{e}{(k)},k\in I,$ be the Fourier transform of $f_{e}(\lambda )$%
:
\begin{equation*}
\hat{f}_{e}{(k)}=\int_{S}e^{ik\lambda }f_{e}(\lambda )\mu (d\lambda ),\quad
k\in I,
\end{equation*}%
where $I=\mathbb{Z}^{d}$ in the torus case and $I=\mathbb{R}^{d}$ in the
case $S=\mathbb{R}^{d}$, respectively. In this last case, we would also need
to assume that $f_{e}\in L_{1}(\mathbb{R}^{d},d\mu )$, for the Fourier
transform to be well defined. Consider the extension of our cumulants:
\begin{eqnarray}
\tilde{J}_{T} &=&\mathrm{Tr}[\prod_{i}^{n}T_{T}(f_{i})]  \notag  \label{matt}
\\
&=&\int_{j_{1},...,j_{n}\in I_{T}}\;\hat{f}_{1}{(j_{2}-j_{1})}\hat{f}_{2}{%
(j_{3}-j_{2})}...\hat{f}_{n}{(j_{1}-j_{n})}\prod_{v=1}^{n}\nu (dj_{v}).
\end{eqnarray}

Replacing the sequences $\hat{f}_{e}(t)$ by their Fourier representations $%
\hat{f}_{e}(t)=\int_{S}f_{e}(\lambda )e^{it\lambda }d\lambda $ in (\ref{matt}%
) yields the following alternative spectral integral representation for
traces of products of Toeplitz matrices or truncated Toeplitz operators%
\footnote{%
Of course, the two expressions $\tilde{J}_{T},J_{T}$ are equal if $f_{e}\in
L_{1},e=1,...,n.$ Note however that the \textquotedblleft spectral
representation\textquotedblright\ (\ref{classp}) is well defined even
without the last condition.}:
\begin{eqnarray}
&&{J}_{T}=\int_{\lambda _{1},...,\lambda _{n}\in S}\;{f}_{1}{(\lambda _{1})}{%
f}_{2}{(\lambda _{2})}...{f}_{n}{(\lambda _{n})}\prod_{e=1}^{n}\Delta
_{T}(\lambda _{e+1}-\lambda _{e})\prod_{e=1}^{n}\mu ({d\lambda _{e}})=
\label{classp} \\
&&\int_{u_{1},...,u_{n-1}\in S}\Big(\int_{\lambda \in S}\;{f}_{1}{(\lambda )}%
{f}_{2}{(\lambda +u_{1})}...{f}_{n}(\lambda +\sum_{1}^{n-1}u_{e})d\lambda %
\Big)\Phi _{T}(u_{1},...,u_{n-1})\prod_{e=1}^{n-1}\mu ({du_{e}}),  \notag
\end{eqnarray}%
where the index $n+1$ is defined to be equal $1$, and where
\begin{equation*}
\Phi _{T}=\Delta _{T}(-\sum_{1}^{n-1}u_{e})\prod_{e=1}^{n-1}\Delta
_{T}(u_{e}),
\end{equation*}%
which, after normalization with the factor $\frac{1}{\left( 2\pi \right)
^{\left( n-1\right) d}T^{d}}$ , yields the \textquotedblleft multivariate Fej%
\'{e}r kernel\textquotedblright\
\begin{equation*}
\Phi _{T}^{\ast }(u_{1},...,u_{n-1})=\frac{1}{\left( 2\pi \right) ^{\left(
n-1\right) d}T^{d}}\Phi _{T}(u_{1},...,u_{n-1}).
\end{equation*}%
Note that the inner integral
\begin{eqnarray}
C(u) &=&C_{(f_{1},...,f_{n})}(u_{1},...,u_{n-1}):=  \notag  \label{Tgcon} \\
&&\int_{\lambda \in S}\;{f}_{1}{(\lambda )}{f}_{2}{(\lambda +u_{1})}...{f}%
_{n}(\lambda +\sum_{1}^{n-1}u_{e})d\lambda ,
\end{eqnarray}%
to be called a graph convolution, is well defined precisely under the
classical H\"{o}lder conditions, when the integrability indicies in (1.2)
satisfy:
\begin{equation}
\begin{cases}
\sum_{e}p_{e}^{-1} & \leq 1\text{ when }S=\mathbb{Z}^{d} \\
\sum_{e}p_{e}^{-1} & =1\text{ when }S=\mathbb{R}^{d},%
\end{cases}
\label{HT}
\end{equation}

The resulting integral
\begin{equation*}
{J}_{T}=\int_{u_{1},...,u_{n-1}\in S}C(u_{1},...,u_{n})\Phi
_{T}(u_{1},...,u_{n-1})\prod_{e=1}^{n-1}\mu ({du_{e}})
\end{equation*}%
is our first example of a \textquotedblleft Fej\'{e}r graph integral" to be
introduced in general in the next section. These are integrals involving
products of Dirichlet kernels $\Delta _{T}$ and functions, applied to linear
combinations, which are related to the vertex-edge incidence structure of a
certain directed graph (in the occurrence, the cyclic graph on the vertices $%
\{1,...,n\}$).

The second expression in the RHS of (\ref{classp}) reveals the asymptotic
behavior of Fej\'er graph integrals, since, as noted, {when $T\rightarrow
\infty $, the multivariate Fej\'er kernel }$\Phi _{T}^{\ast }$ \textbf{%
convergences weakly to a $\delta $ measure}.

The kernel lemma \ref{lwc} will imply that

\begin{equation*}
T^{-d}{J}_{T}\rightarrow \left( 2\pi \right) ^{\left( n-1\right)
d}C(0,...,0)=\int_{\lambda \in S}\;{f}_{1}{(\lambda )}{f}_{2}{(\lambda )}...{%
f}_{n}(\lambda )\mu (d\lambda )\text{ \ \ }
\end{equation*}
provided that we check that the function $C(u_1, ...,u_{n-1})$ is bounded
and continuous. This is indeed true, as stated in the next result:

\begin{lmm}
\label{l:contconv} The ``graph convolution" function $C(u_1, ...,u_{n-1})$
defined in (\ref{Tgcon}) is bounded and continuous if (\ref{e:Lp1}) holds
with integrability indices satisfying the condition (\ref{HT}).
\end{lmm}

\textbf{Proof:} Note that the function $C:S^{n-1}\rightarrow \mathbb{R}$ is
a composition
\begin{equation*}
C(u_{1},...,u_{n-1})=\mathcal{C}%
_{(T_{1}(u_{1},...,u_{n-1}),...,T_{n}(u_{1},...,u_{n-1}))}
\end{equation*}%
of the functional
\begin{equation*}
\mathcal{C}_{(f_{1},...,f_{n})}:\prod_{e=1}^{n}{\mathbf{L}}%
_{p_{e}}\rightarrow \mathbb{R}.
\end{equation*}%
defined by \footnote{%
Note again this is well defined precisely under the classical H\"{o}lder
conditions.}
\begin{equation*}
\mathcal{C}_{(f_{1},...,f_{n})}:=\int_{\lambda \in
S}\prod_{e=1}^{n}f_{e}(\lambda )\mu (d\lambda )
\end{equation*}%
with the continuous functionals
\begin{equation*}
T_{e}(u_{1},...,u_{n-1}):S^{n-1}\rightarrow {{\mathbf{L}}}_{p_{e}}
\end{equation*}%
defined by $T_{e}(u_{1},...,u_{n-1})=f_{e}(\cdot
+\sum_{v=1}^{e-1}u_{v}),e=1,...,n$.

Indeed, the continuity of the functionals $T_e$ is clear when $f_e$ is a
function which is continuous, bounded and of bounded support, and this
continues to be true for functions $f_e \in {{\mathbf{L}}}_{p_e}$, since
these can be approximated in the $L_{p_e}$ sense. In conclusion, under the
H\"older assumptions, the continuity of the functional $C(u_1,...,u_{n-1})$
will follow automatically from that of $\mathcal{C}_{(f_1,...,f_n)}$.

Finally, under the \textquotedblleft H\"{o}lder conditions" (\ref{HT}), the
continuity as well as boundedness of the function $\mathcal{C}%
_{(f_{1},...,f_{E})}$ follow from its multi-linearity and from the H\"{o}%
lder inequality:
\begin{equation*}
|\mathcal{C}_{(f_{1},...,f_{n})}|\leq \prod_{e=1}^{n}\Vert {f_{e}}\Vert
_{p_{e}}\ .
\end{equation*}

\begin{thrm}
\label{t:Toep} Let $f_{e}\in \mathbf{L}_{p_{e}},e=1,...,n$ where $1\leq
p_{e}\leq \infty ,e=1,...,n$, and let $J_{T},\tilde{J}_{T}$ be defined by (%
\ref{classp}), (\ref{matt}) respectively. Then, it follows that:

\begin{enumerate}
\item If the integrability indices satisfy the H\"{o}lder conditions
\begin{equation*}
\begin{cases}
\sum_{e}p_{e}^{-1} & \leq 1\text{ when }S=\mathbb{Z}^{d} \\
\sum_{e}p_{e}^{-1} & =1\text{ when }S=\mathbb{R}^{d},%
\end{cases}%
\end{equation*}%
then
\begin{equation}
\lim_{T\rightarrow \infty }T^{-d}J_{T}=\int_{\lambda \in
S}\prod_{e=1}^{n}f_{e}(\lambda )\mu (d\lambda )  \label{e:T}
\end{equation}

\item If $\alpha :=\sum_{e}p_{e}^{-1}>1,$ then it holds that:
\begin{equation*}
J_{T}=o(T^{\alpha d}).
\end{equation*}

\item If in the continuous case it holds in addition that if $f_e \in
\mathbf{L}_{p_e} \cap L_1, e=1,..., n$, then Fourier coefficients $\hat{f}%
_e(k)$ may be defined, and the previous results hold also for $\tilde{J}_T=%
\mathrm{Tr}[\prod_{e=1}^{n}T_{T}(f_{e})]$.
\end{enumerate}
\end{thrm}

\textbf{Notes:} 1) In the second case, the exact exponent of magnitude is
unknown, except for the upper bound $\alpha \,d$.

2) The result above is a refinement of a limit theorem of Grenander and Szeg%
\"{o} \cite{grenander:szego:1958} concerning traces of products of truncated
Toeplitz operators. Under the current strengthened integrability conditions,
it was obtained when $d=1$ in the discrete case in \cite{a} and in the
continuous case in \cite{Gin}; as we show below, the result holds in fact in
any dimension $d$ (after modifying the necessary integrability condition for
the asymptotic variance, in accordance to the H\"older-Young-Brascamp-Lieb
inequality).

Following the proof of \cite{AB}, we see that part 1 of Theorem \ref{t:Toep}
follows from Lemmas \ref{lwc}, \ref{l:contconv} above, yielding the
convergence of the normalized variance:
\begin{equation*}
\lim_{T\rightarrow \infty }\frac{\chi _{2,T}}{T^{d}}=J(0,...,0)=\int_{%
\lambda \in S}\;{f}_{1}{(\lambda )}{f}_{2}{(\lambda )}...{f}_{k}{(\lambda )}%
\mu (d\lambda ).
\end{equation*}%
Part 2 (see the proof of Theorem 3.1(c) below) implies that for $k\geq 3,$
the cumulants satisfy
\begin{equation*}
\lim_{T\rightarrow \infty }\frac{\chi _{k,T}}{T^{dk/2}}=0,
\end{equation*}%
and implies asymptotic normality. We arrive thus at the following
multidimensional generalization of the results of Avram \cite{Avram1988} and
Ginovian \cite{Gin}.

\begin{thrm}
\label{t:quadr} Consider the quadratic functional

\begin{equation*}
Q_{T}=Q_{T}^{(1,1)}=\int_{t,s\in I_{T}}\left[ X_{t}X_{s}-{\mathbb{E}}%
X_{t}\,X_{s}\right] \hat{b}(t-s)\nu (ds)\nu (dt),
\end{equation*}
where $X_{t},\quad t\in {\ I}$ is a Gaussian random field with spectral
density $f(\lambda )\in \mathbf{L}_{p}$. Assume the generating function $b$
of the quadratic functional is such that $b(\lambda )\in \mathbf{L}_{q},$
that:
\begin{equation*}
\frac{1}{p}+\frac{1}{q}\leq \frac{1}{2},
\end{equation*}
and that in the continuous case we have also $b(\lambda ), f(\lambda) \in
L_{1}.$

Then, the central limit theorem holds:
\begin{equation*}
\lim T^{-d/2}Q_{T}\rightarrow N(0,\sigma ^{2}),\ \ T\rightarrow \infty ,
\end{equation*}

where
\begin{equation}  \label{sG}
\sigma ^{2}:=2(2\pi )^{d}\int_{S}b^{2}(\lambda )f^{2}(\lambda )d\lambda .
\end{equation}
\end{thrm}

We present now one more result for Gaussian fields, which is related in the {%
discrete} case to the classical result of Breuer and Major \cite{BM}, and in
the continuous case to the result of Ivanov and Leonenko \cite%
{IvanovLeonenko1989}. We note that these authors worked under time-domain
assumptions, however, reasoning in the spectral domain with the methodology
of \cite{a} and the present paper -- see Example \ref{exagraph} --
immediately lead to the following result.

\begin{thrm}
\label{qs} Let $X_{t}$, $t\in {I}$, be a Gaussian random field with spectral
density $f(\lambda )\in \mathbf{L}_{p}$. Let $S_{T}=\int_{t\in
I_{T}}P_{l}(X_{t})\nu (dt)$, where $P_{l}(X_{t})$ are univariate Appell
(Hermite) polynomials and $l\geq 2$. Assume that:
\begin{equation}  \label{hys}
z:=p^{-1}\leq 1-\frac{1}{l}.
\end{equation}

Then,
\begin{eqnarray}  \label{varci}
\sigma^2:= f^{(*,l)}(0)=\int_{y_1, ...,y_{l-1} \in S} f(y_1)
f(y_2-y_1)...f(y_{l-1}-y_{l-2}) f(-y_{l-1}) \prod_{i=1}^{l-1} dy_i < \infty.
\end{eqnarray}

If, moreover, $\sigma \neq 0$, then the central limit theorem holds
\begin{equation*}
\lim T^{-d/2}S_{T}\rightarrow N(0,\sigma ^{2}).
\end{equation*}
\end{thrm}

\textbf{Note:} The difference between the integral representations of the
variances (\ref{sG}) and (\ref{varci}) will be explained via graph theory in
the next section.


\section{Fej\'er graph/matroid integrals and graph/matroid convolutions\label%
{s:kgi}}

In this section, we introduce a unifying graph-theoretical framework for
problems similar to those of the previous section.

\begin{dfntn}
\label{ss:delta} \textit{Let $G=({\mathcal{V}, \mathcal{E}})$ denote a
directed graph with $V$ vertices, $E $ edges, a basis of $C$ independent
cycles \footnote{%
a basis of cycles is a set of cycles, none of which may be obtained via
addition modulo $2$ of other cycles, after ignoring the orientation} and $%
co(G)$ components. \textbf{The incidence matrix} $M=\{M_{v,e}\}_{v \in
\mathcal{V}, e \in \mathcal{E}}$ of the graph is the $V \times E$ matrix
with entries $[v,e]=\pm 1 $ if the vertex $v$ is the end/start point of the
edge $e$, and $0$ otherwise.}

A \textbf{circuit matrix} $M^*$ is a $C \times E$ matrix whose rows are
obtained by assigning arbitrary orientations to a basis of circuits (cycles)
$c=1,...,C$ of the graph, and by writing each edge as a sum of $\pm$ the
circuits it is included in, with the $\pm$ sign indicating a coincidence or
opposition to the orientation of the cycle \label{dinc}
\end{dfntn}

Besides the graph framework, we will hint also to possible matroid
generalizations. To clarify this point, let us start by quoting Tutte:
\textquotedblleft it is probably true that any theorem about graphs
expressible in terms of edges and circuits exemplifies a more general result
about vector matroids".

Let us recall briefly that matroids are a concept which formalizes the
properties of the \textquotedblleft rank function" $r(A)$ obtained by
considering the rank of an arbitrary set of columns $A$ in a given arbitrary
matrix $M$. More precisely, a matroid is a pair ${\mathcal{E}},r:2^{\mathcal{%
E}}\rightarrow \mathbb{N}$ of a set ${\mathcal{E}}$ and a \textquotedblleft
rank like function" $r(A)$ defined on the subsets of ${\mathcal{E}}$.
Matroids may also be defined in equivalent ways via their independent sets,
via their bases (maximal independent sets), via their circuits (minimal
dependent sets), via their spanning sets (sets containing a basis), or via
their flats (sets which may not be augmented without increasing the rank).
For precise definitions and for excellent expositions on graphs and
matroids, see \cite{Ox}, \cite{Ox1} or \cite{W}.

The most familiar matroids, called \emph{vectorial matroids}, are defined by
the set ${\mathcal{E}}$ of columns of a matrix and by the rank function $%
r(A) $ which gives the rank of any set of columns $A$ (matrices with the
same rank function yield the same matroid).

Some useful facts from matroid theory are the fact that to each matroid $M$
one may associate a \textbf{dual matroid}, with rank function
\begin{equation*}
r^{\ast }(A)=|A|-r(M)+r({\mathcal{E}}-A).
\end{equation*}

For vectorial matroids, the dual matroid is also vectorial, associated to
any matrix whose rows span the space orthogonal to the rows of $M$.
Furthermore, in the case of graphic matroids, the dual matroid is associated
to the circuit matrix.

Tutte's \textquotedblleft conjecture" holds true in our case: a matroid Szeg%
\"{o}-type limit theorem was already given in \cite{a}. However, for
simplicity, we will restrict ourselves here to the particular case of
\textbf{graphic matroids} associated to the incidence matrix $M$ of an
oriented graph. In this case, the proofs are more intuitive, due to the fact
that the algebraic dependence structures translate into graph-theoretic
concepts, like circuits corresponding to cycles, etc.

From here on, we will restrict ourselves to the graphic case, i.e. to the
case when our dependence matrix is the incidence matrix of a directed graph.

Let $(S,d\mu )$ denote either $\mathbb{R}^{d}$ with Lebesgue measure, or the
torus $[-\pi ,\pi ]^{d}$ with normalized Lebesgue measure, and let ${f}%
_{e}(\lambda ):S\rightarrow \mathbb{R},\ e=1,\ldots ,E$ denote a set of
functions associated to the columns of $M$, which satisfy integrability
conditions
\begin{equation}
f_{e}\in {\mathbf{L}}_{p_{e}}(S,d\mu ),\ 1\leq p_{e}\leq \infty .
\label{e:Lp}
\end{equation}%
%
%
%
%
%
%
Let $\hat{f}_{e}{(k)},k\in I,$ the Fourier transform of $f_{e}(\lambda )$:
\begin{equation*}
\hat{f}_{e}{(k)}=\int_{S}e^{ik\lambda }f_{e}(\lambda )\mu (d\lambda ),\quad
k\in I,
\end{equation*}%
where $I=\mathbb{Z}^{d}$ in the torus case and $I=\mathbb{R}^{d}$ in the
case $S=\mathbb{R}^{d}$, respectively. In this last case, we would also need
to assume that $f_{e}\in L_{1}(\mathbb{R}^{d},d\mu )$, for the Fourier
transform to be well defined. However, all our analytic results concern the
spectral domain, and hence this assumption will not be necessary.

Our object of interest, in its \textquotedblleft time domain
representation", is:
\begin{eqnarray}
\tilde{J}_{T} &=&\tilde{J}_{T}(M,f_{e},e=1,...,E)  \notag \\
&=&\int_{j_{1},...,j_{V}\in I_{T}}\;\hat{f}_{1}{(i_{1})}\hat{f}_{2}{(i_{2})}%
...\hat{f}_{E}{(i_{E})}\prod_{v=1}^{V}\nu (dj_{v}),  \label{matF}
\end{eqnarray}%
where $i=(i_{1},...,i_{E})=(j_{1},...,j_{V})M=jM$, where $\nu (dj_{v})$
stands for Lebesgue measure and counting measure, respectively, and where in
the torus case the linear combinations are computed modulo $[{-\pi },{\ \pi }%
]^{d}$, so that the linear map $jM:S^{V}\rightarrow S^{E}$ is well defined.

\textbf{Note.} To keep the transparent analogy with the case, when $d=1$, we
make the following convention concerning notations. Here and in what follows
let us treat a product of a vector, whose components are $d$-dimensional,
and a matrix (or another vector) with scalar components in a specific sense:
we will still perform multiplication component-wise according to the usual
rule, and as a result we obtain a vector, whose components are $d$%
-dimensional again (or, correspondingly, just $d$-dimensional vector).

A \textbf{Fej\'er graph integral} is the expression obtained by replacing
the sequences $\hat{f}_{e}(t)$ in (\ref{matF}) by their Fourier
representations $\hat{f}_{e}(t)=\int_{S}f_{e}(\lambda )e^{it\lambda }
\mu(d\lambda) $: under the assumption $f_e \in L_{1} (S, d \mu )$, an easy
computation (see \cite{AB}, Lemma 1) shows that (\ref{matF}) may be written
also as the integral (\ref{matF2}) below. We introduce however a more
general concept.

\begin{dfntn}
\label{d:delta} \textit{Let $(S, d \mu)$ denote either $\mathbb{R}^d$ with
Lebesgue mesure, or the torus $[-\pi,\pi]^d$ with normalized Lebesgue
mesure. Let $M$ be a matrix of dimensions $V \times E$, with arbitrary
coefficients in the first case and with integer coefficients in the second
case. Let ${f}_e(\lambda): S \to \mathbb{R},\ e=1,\ldots,E$ denote a set of
functions associated to the columns of $M$. Suppose these functions satisfy
integrability conditions
\begin{equation}  \label{e:Lp}
f_e \in {\mathbf{L}}_{p_e} (S, d \mu ),\ 1\leq p_e\leq \infty,
\end{equation}
}

\textit{A Fej\'er matroid integral is defined by the following ``spectral
representation":
\begin{eqnarray}
J_T&=& J_T(M, f_e, e=1,..., E)  \label{matF2} \\
&=&\int_{\lambda_1, ...,\lambda_E \in S} \; {f}_{1}{(\lambda_1)} {f}_{2}{%
(\lambda_2)}...{f}_{E}{(\lambda_E)} \prod_{v=1}^V \Delta_T( u_v)
\prod_{e=1}^E \mu({d \lambda_e})  \notag
\end{eqnarray}
where $\Delta_T(u)$ is a kernel defined by (\ref{ker}), where $%
(u_1,...,u_V)^{\prime}= M (\lambda_1, ..., \lambda_E)^{\prime}$, and where
in the torus case the linear combinations are computed modulo $[{- \pi}, {\
\pi}]^d$. }

\textit{A Fej\'er matroid integral will be called a Fej\'er graph integral
for graphic matroids associated to the incidence matrix $M$ of a directed
graph $G$. In this case, the functions and kernels in (\ref{matF2}) are
associated respectively to the edges and vertices of the graph.}
\end{dfntn}


\textbf{The cycle graph/Toeplitz example} Consider the particular case of a
cyclic graph with $n$ edges. In this case, the matrix $M$ with $n$ columns
and rows, is:
\begin{equation*}
M=\left(
\begin{array}{ccccccc}
-1 & 0 & 0 & . & ... & 0 & 1 \\
1 & -1 & 0 & 0 & ... & 0 & 0 \\
0 & 1 & -1 & \ddots & ... & \vdots & \vdots \\
\vdots &  & \ddots & \ddots &  &  &  \\
0 & 0 &  & \ddots &  & -1 & 0 \\
0 & 0 & . & . & . & 1 & -1%
\end{array}%
\right)
\end{equation*}%
and its Fej\'{e}r graph integral is given by:
\begin{equation*}
J_{T}=\int_{\lambda _{1},...,\lambda _{n}\in S}\;{f}_{1}{(\lambda _{1})}{f}%
_{2}{(\lambda _{2})}...{f}_{n}{(\lambda _{n})}\prod_{v=1}^{n}\Delta
_{T}(\lambda _{v}-\lambda _{v+1})\prod_{e=1}^{n}\mu ({d\lambda _{e}}).
\end{equation*}

\textbf{Note.} For analytical results concerning only Fej\'er matroid
integrals as defined by (\ref{matF2}), the condition $f_e \in L_{1} (S, d
\mu )$ is unnecessary.

\section{Limit theory for Fej\'er graph integrals \label{s:GS}}

The main points of the limit theory for Fej\'er graph integrals, to be
presented now, are that:

\begin{enumerate}
\item Under certain H\"{o}lder-Young-Brascamp-Lieb conditions necessary to
ensure the existence of the limiting integral, the following convergence
holds as $T\rightarrow \infty $:
\begin{equation}
\boxed{ T^{-d} \; J_T(M,f_1,...,f_E) \to \int_{S^C} {f}_{1}{(\lambda_1)}
{f}_{2}{(\lambda_2)}...{f}_{E}{(\lambda_E)} \prod_{c=1}^C \mu(d y_c) },
\label{e:GS}
\end{equation}%
where $(\lambda _{1},...\lambda _{E})=(y_{1},...,y_{C})M^{\ast }$ (with
every $\lambda _{e}$ reduced modulo $[-\pi ,\pi ]^{d}$ in the torus case), $%
M^{\ast }$ being any matrix whose rows span the space orthogonal to the rows
of $M$, and $C$ being the rank of $M^{\ast }$. Informally, the kernels
disappear in the limit, giving rise to the \textquotedblleft dual matroid" $%
M^{\ast }$.

\item When the H\"{o}lder-Young-Brascamp-Lieb conditions do not hold, then,
cf. part c) of the theorem below, the normalization defined in part a) will
lead to a zero limit.
\end{enumerate}

\begin{thrm}
\label{Sz} Suppose that $f_{e}\in {L}_{p_{e}}(d\mu )$ for part a), and $%
f_{e}\in {{\mathbf{L}}}_{p_{e}}(d\mu )$ for parts b),c) and set $%
z=(p_{1}^{-1},...,p_{E}^{-1})$.

Let $J_T=J_T(M, f_1,...,f_E)$ denote a Fej\'er matroid integral and let $%
r(A), r^*(A)$ denote respectively the ranks of a set of columns in $M$ and
in the dual matroid $M^*$.

Suppose that for every row $l$ of the matrix $M$, one has $r(M)=r(M_l)$,
where $M_l$ is the matrix with the row $l$ removed. Then:

\begin{enumerate}
\item[a)]
\begin{eqnarray}  \label{e:bound}
\boxed{J_T(M, f_1,...,f_E) \leq c_M T^{d \; \a_M( z)}}
\end{eqnarray}
where $c_M$ is a constant independent of $z$ and

\begin{equation*}
\begin{cases}
\text{in the discrete case } & \alpha _M( z) \text{ is given by } (\ref%
{borne}) \\
\text{in the continuous case} & \alpha _M( z)= \alpha ^c_M(z)=co(M)+ (\sum_e
z_e -C)_+%
\end{cases}%
\end{equation*}

\item[b)] If $\alpha _M( z)=V-r(M)=co(M)$, or, equivalently,
\begin{eqnarray}  \label{pcc}
&\sum_{j \in A} z_j \leq r^*(A),\; \ \forall A &\text{ in the discrete case}
\\
&\sum_e z_e \leq C \; &\text{ in the continuous case}  \notag
\end{eqnarray}
then

\begin{equation}
\lim_{T\rightarrow \infty }\frac{J_{T}(M)}{T^{d\;co(M)}}=k_{M}\;\mathcal{J}%
(M^{\ast },f_{1},...,f_{E}),  \label{t:GS}
\end{equation}%
where
\begin{equation}
\!\!\!\!\!\boxed{ \sJ(M^*,f_1,...,f_E) = \int_{S^C} {f}_{1}{(\lambda_1)}
{f}_{2}{(\lambda_2)}...{f}_{E}{(\lambda_E)} \prod_{c=1}^C \mu(d y_c) }
\label{limint}
\end{equation}%
and where $(\lambda _{1},...\lambda _{E})=(y_{1},...,y_{C})M^{\ast }$ (with
every $\lambda _{e}$ reduced modulo $[-\pi ,\pi ]^{d}$ in the discrete
case), and $C$ denotes the rank of the dual matroid $M^{\ast }$.

\item[c)] If a strict inequality $\alpha _{M}(z)>co(M)$ holds, then the
inequality (\ref{e:bound}) of Theorem \ref{Sz} a) may be strengthened to:
\begin{equation*}
J_{T}(M)=o(T^{d\;\alpha _{M}(z)})
\end{equation*}
\end{enumerate}
\end{thrm}

\noindent \textbf{Remark:} The results of this theorem, that is, the
expression of $\alpha _M( z)$ and the limit integral $\mathcal{J}(M^*)=%
\mathcal{J}(M^*,f_1,...,f_E)$, as well as the convergence conditions of
integrals depend on $M, M^*$ only via the two equivalent rank functions $%
r(A), r^*(A)$, i.e. only via the matroid dependence structure between the
columns, and not on the chosen representing matrices.


\textbf{Proof:} The proof of part b) of Theorem \ref{Sz} is essentially
identical with that given in \cite{AB}, up to the modification of the
integrability conditions and the appearance of the extra constant $k_M.$ For
completeness, we sketch now this proof, for a connected graph (w.l.o.g.).

Note first that in a connected graph there are only $V-1$ independent rows
of the incidence matrix $M$ (or independent variables $u_j$), since the sum
of all the rows is $0$ (equivalently, $u_V =-\sum_{v=1}^{V-1} u_v$). Thus, $%
r(M)=V-1$, $co(M)=1$, and the order of magnitude appearing in the
normalization is just $T^d$.

The main idea behind the proof of Theorem \ref{Sz} b) are a change of
variables and applying the continuity of graph convolutions:

\begin{enumerate}
\item \textbf{Change of variables.} Fix a basis $y_{1},...,y_{C}$ in the
complement of the space generated by the $u_{v}$'s, $v=1,\ldots ,V$, switch
to the variables $u_{1},...,u_{V-1},y_{1},...,y_{C}$ and integrate in (\ref%
{matF2}) first over the variables $y_{c}$'s, $c=1,\ldots ,C$. This is more
convenient in the graphic case, since, after fixing an arbitrary spanning
tree ${\mathcal{T}}$ in the graph, the complementary set of edges ${\mathcal{%
T}}^{c}$ furnishes a maximal set of independent cycles (with cardinality $C$%
). Assume w.l.o.g. that in the list $(\lambda _{1},...,\lambda _{E})$, the
edges in ${\mathcal{T}}^{c}$ are listed first, namely $(\lambda _{e},e\in {%
\mathcal{T}}^{c})=(\lambda _{1},...,\lambda _{C})$. We make the change of
variables $y_{1}=\lambda _{1},...,y_{C}=\lambda _{C}$, and $%
(u_{1},...,u_{V-1})^{\prime }=\tilde{M}(\lambda _{1},...,\lambda
_{E})^{\prime }$, where $\tilde{M}$ denotes the first $V-1$ rows of the
incidence matrix $M$. Thus,
\begin{equation*}
(y_{1},...,y_{C},u_{1},...,u_{V-1})^{\prime }=%
\begin{pmatrix}
I_{C} & 0 \\
\tilde{M}_{C} & \tilde{M}_{V}%
\end{pmatrix}%
\;(\lambda _{1},...,\lambda _{E})^{\prime }
\end{equation*}%
where the first rows are given by an identity matrix $I_{C}$ completed by
zeroes and where $\tilde{M}_{C},\tilde{M}_{V}$ denote the first $C$ columns/
next $V-1-C$ columns of the matrix $\tilde{M}$.

Inverting the transformation above yields:
\begin{equation}  \label{e:inv}
(\lambda_1, ..., \lambda_E)=(y_1, ..., y_C, u_1, ...,u_{V-1}) \;
\begin{pmatrix}
I_C & - \tilde{M}_{C}^{\prime}\, \tilde{M}_{V}^{-1} \\
0 & \tilde{M}_{V}^{-1}%
\end{pmatrix}%
= (y_1, ..., y_C, u_1, ...,u_{V-1}) \;
\begin{pmatrix}
{M}^* \\ \hline
N%
\end{pmatrix}%
,
\end{equation}
that is, it turns out that the first rows of the inverse matrix are
precisely the dual matroid $M^*$.

\begin{dfntn}
The function
\begin{equation}  \label{e:mconv}
h_{M^*,N}( u_1, ...,u_{r(M)}) =\int_{y_1, ..., y_C \in S} {f}_{1}{(\lambda_1)%
} {f}_{2}{(\lambda_2)}...{f}_{E}{(\lambda_E)} \prod_{c=1}^C d \mu(y_c)
\end{equation}
where $\lambda_e$ are represented as linear combinations of $y_1, ..., y_C,
u_1, ...,u_{V-1}$ via the linear transformation (\ref{e:inv}) will be called
a \textbf{matroid/graph convolution} depending on whether the matroid is
graphic or not.
\end{dfntn}

The change to the variables $y_{1},...,y_{C},u_{1},...,u_{V-1}$ and
integration over $y_{1},...,y_{C}$ transforms the Fej\'{e}r graph integral
into the following integral of the product of a \textquotedblleft graph
convolution" and a Fej\'{e}r kernel:
\begin{equation*}
J_{T}(M)=\int_{u_{1},...,u_{V-1}\in S}h_{M^{\ast
},N}(u_{1},...,u_{V-1})\prod_{v=1}^{V}\Delta
_{T}(u_{v})\;\prod_{v=1}^{V-1}d\mu (u_{v}).
\end{equation*}

Recalling that the Fej\'{e}r kernel converges under appropriate conditions
to Lebesgue measure on the set $u_{1}=...=u_{V-1}=0$, we find, just as in
the cycle case, that part b) of Theorem \ref{Sz} will be established once
the convergence of the kernels and the continuity of the graph convolutions $%
h(u_{1},...,u_{r(M)})$ in the variables $(u_{1},...,u_{r(M)})$ is
established.

\item
\begin{lmm}
\textbf{The continuity of graph convolutions.} \label{lbcg} The
\textquotedblleft graph convolution" function $C(u_{1},...,u_{n-1})$ defined
in (\ref{Tgcon}) is bounded and continuous if (\ref{e:Lp1}) holds with
integrability indices satisfying the power counting condition (1.6).
\end{lmm}

The proof is essentially the same as in the cycle case. Note that the
function $h:\mathbb{R}^{V-1}\rightarrow \mathbb{R}$ is a composition
\begin{equation*}
h_{M^{\ast },N}(u_{1},...,u_{V-1})=J(M^{\ast },T_{1}(f_{1}),...,T_{E}(f_{E}))
\end{equation*}%
of the continuous functionals
\begin{equation*}
T_{e}(u_{1},...,u_{V-1}):\mathbb{R}^{V-1}\rightarrow {{\mathbf{L}}}_{p_{e}}
\end{equation*}%
and of the functional
\begin{equation*}
J(M^{\ast },f_{1},...,f_{E}):\prod_{e=1}^{E}{{\mathbf{L}}}%
_{p_{e}}\rightarrow \mathbb{R}.
\end{equation*}%
The functional $T_{e}$ is defined by $T_{e}(u_{1},...,u_{V-1})=f_{e}(\cdot
+\sum_{v}u_{v}N_{v,e})$, where the $N_{v,e}$ are the components of the
matrix $N$ in (\ref{e:inv}). The functionals $T_{e}$ are clearly continuous
when $f_{e}$ is a continuous function, and this continues to be true for
functions $f_{e}\in {{\mathbf{L}}}_{p_{e}}$, since these can be approximated
in the $L_{p_{e}}$ sense by continuous functions. Thus, under our
assumptions, the continuity of the functional $h_{M^{\ast
},N}(u_{1},...,u_{V-1})$ follows automatically from that of $J(M^{\ast
},f_{1},...,f_{E})$.

Finally, under the ``power counting conditions" (\ref{pcc}), the continuity
of the function $J(M^*,f_1,...,f_E)$ follows from the
H\"older-Brascamp-Lieb-Barthe inequality:
\begin{equation*}
|J(M^*,f_1,...,f_E)| \leq \prod_{e=1}^E \Vert {f_e} \Vert_{p_e}
\end{equation*}
(see Theorem \ref{t:HYB}).

In conclusion, the convergence of the Fej\'er kernels to a $\delta $ measure
implies the convergence of the scaled Fej\'er graph integral
\begin{equation*}
J_T(M, f_e, e=1...,E) \quad \text{to} \quad \mathcal{J}(M^*, f_e, e=1...,E),
\end{equation*}
establishing Part b) of the theorem.
\end{enumerate}

The proof of parts a), c) are postponed to section 3.2.

\noindent \textbf{Remarks:} 1) In the spatial statistics papers (\cite{ALS},
\cite{als}), the continuity of the graph convolutions $h_{M^*,N}(u_1,
...,u_{V-1})$ was assumed to hold, and indeed checking whether this
assumption may be relaxed was one of the outstanding difficulties for the
spatial extension.

2) It is not difficult to extend this approach to the case of several
components and then to the matroid setup. In the first case, one would need
to choose independent cycle and vertex variables $y_{1},...,y_{r(M^{\ast })}$
and $u_{1},...,u_{r(M)}$, note the block structure of the matrices, with
each block corresponding to a graph component, use the fact that for graphs
with several components, the rank of the graphic matroid is $r(M)=V-co(G)$
and finally Euler's relation $E-V~=C-co(G)$, which ensures that
\begin{equation*}
E=\big(V-co(G)\big)+C=r(M)+r(M^{\ast }).
\end{equation*}

3) An important feature of the discrete case is that the limiting result
(relation (\ref{t:GS})) when $f_e$ are complex exponentials is
straightforward, implying therefore immediately theorem \ref{Sz} in this
case, by the multilinearity of $\mathcal{J}(M^*,f_1,...,f_E) $ and $T^{-d}
J(M,f_1,...,f_E) $ and by Lemma \ref{l:exp} below.

\begin{lmm}
\label{l:exp} For any matrix with $r(M)=V-1$, any set of integers $b=(b_e,
e=1, ...,E)$, and any functions $f_e(\lambda_e)=e^{i \lambda_e b_e}, e=1,
...,E$, theorem \ref{Sz} holds, i.e.:
\begin{equation*}
\lim_{T\rightarrow \infty }\frac{\int_{S^E} e^{ i <b,\lambda >}
\prod_{v=1}^V \Delta_T(u_v) \prod_{e=1}^E \mu(d \lambda_e)}{T^{d }} = k_M \;
\int_{S^C} e^{ i <\lambda ,b>} \prod_{c=1}^C \mu(d y_c)
\end{equation*}
\end{lmm}

\textbf{Proof:} Aside from the constant $k_M$, the stated RHS (limiting
value) above is:
\begin{eqnarray*}
\int_{S^C} e^{ i <\lambda ,b>} \prod_{c=1}^C \mu(d y_c)=\int_{S^C} e^{ i <y,
M^* b>} \prod_{c=1}^C \mu(d y_c)=1_{M^* b=0}=1_{ b \in R(M)}
\end{eqnarray*}
where $R(M)$ denotes the subspace generated by the rows of $M$.

Now the LHS in (\ref{t:GS}), before scaling, is:
\begin{eqnarray*}
\int_{S^E} e^{ i <b,\lambda >} \int_{I_T^V} e^{ - i <s M, \lambda >}
\prod_{v=1}^V \nu(d s_v) \prod_{e=1}^E \mu(d \lambda_e) =\int_{I_T^V}
\prod_{v=1}^V \nu(d s_v) (\int_{S^E} e^{ i <b- s M,\lambda>} \prod_{e=1}^E
\mu(d \lambda_e)) \\
=\int_{I_T^V} \prod_{v=1}^V \nu(d s_v) 1_{ b-s M=0}= \nu(s: s M =b, s \in
I_T)
\end{eqnarray*}

Now for a matrix with integer entries it holds that:
\begin{equation*}
E(T):=\nu (s:sM=b,s\in I_{T})
\end{equation*}%
is either
\begin{equation*}
E(T)%
\begin{cases}
=0 & \text{ if }b\notin R(M), \\
\sim T^{D}\mu (I_{1}\cap ker(M)) & \text{ if }b\in R(M),%
\end{cases}%
\end{equation*}%
the second statement being tantamount to the definition of the Lebesgue
measure. Thus, the result holds with $k_{M}=\mu (I_{1}\cap ker(M))$. See for
more details \cite{a}.

Note also that $E(T)$ is an ``Ehrhart quasi polynomial", whose next
coefficients are related to other geometric characteristics of $I_1$, which
should allow developing correction terms to Theorem \ref{Sz}.

\textbf{Note:} The lemma above may be interpreted as saying that the
measures on $(S)^{E}$ given by the \textquotedblleft multiple Fej\'{e}r
kernels"
\begin{equation*}
T^{-d}\Delta _{T}(-\sum_{v=1}^{V-1}\lambda _{v})\prod_{v=1}^{V-1}(\Delta
_{T}(\lambda _{v})\nu (d\lambda _{\nu }))
\end{equation*}%
converge weakly as $T\rightarrow \infty $ to the uniform measure on the
subspace $\lambda =yM^{\ast }$ (since the Fourier coefficients converge).

\subsection{The upper bound for the order of magnitude of Fej\'er matroid
integrals \label{s:GS1}}

We turn now to the \textquotedblleft upper bound exponent" $\alpha _{M}(z)$
for the order of magnitude of Fej\'{e}r matroid integrals (useful when it is
not precisely $d$). In the discrete case \cite{a}, the exponent $\alpha
_{M}(z)$ of this upper bound turns out to be $d$ times the solution of a
graph optimization problem:
\begin{equation}
\boxed{\a_M( z) = co(M) + \max_{A \subset {1,...,E} } [\sum_{j \in A} z_j -
r^*(A)] }  \label{borne}
\end{equation}%
or equivalently,
\begin{equation}
\boxed{ \a_M( z)= \max_{A \subset {1,...,E}} [co(M-A) - \sum_{j \in A}(1-
z_j)] }  \label{borne2}
\end{equation}%
where $co(M-A)$ represents the number of remaining components, after the
edges in $A$ have been removed, and, for a general Fej\'{e}r graph integral,
we define
\begin{equation}
co(M-A)=V-r(M-A).  \label{e:co}
\end{equation}

Note that for connected graphs and under the power counting conditions $%
\sum_{j \in A} z_j \leq r^*(A)$, this exponent reduces to $d$, as in Theorem %
\ref{t:Toep}.

We will call the problem (\ref{borne2}) a \textbf{graph breaking problem:
find a set of edges whose removal maximizes the difference between the
number of remaining components and} $\sum_{j \in A} (1-z_j) $, or,
equivalently, the difference between $\sum_{j \in A}z_j$ and the dual rank $%
r^*(A)$.

\textbf{Note:} In the following examples an important role is played by the
\textbf{``maximal breaking" $A=M$} and \textbf{``no breaking $A=\emptyset$"}
sets $A$, which yield often the solution of the optimal breaking problem. It
is useful to introduce therefore the lower bound:

\begin{eqnarray}
&&\alpha^c_M( z)= \max_{A \in \{\mathcal{E}, \emptyset\}} [co(M-A) - \sum_{j
\in A}(1- z_j)]= \max\{co(M),\sum_e (z_e-1) + V\}  \notag \\
&&=\max\{co(M),co(M) + \sum_e z_e -C\} =co(M)+ (\sum_e z_e -C)_+ \ ,
\label{ub}
\end{eqnarray}
where the equality before the last holds by Euler's relation $C=(E-V)_+$ ($C$
denotes the number of cycles, and, more generally, the rank of the dual
matroid $M^*$).

Note that in the case of a cycle graph of size $m$, this reduces to
\begin{equation*}
\alpha^c_m( z)=\max\{1,\sum_{e=1}^m z_e\},
\end{equation*}
as stated in Theorem \ref{t:Toep}, and that the expression (\ref{ub}) turns
out to yield the upper bound exponent in the continuous case.

\textbf{Note:} Theorem \ref{Sz} c) may be used for establishing convergence
to $0$ of higher order cumulants, whenever these may be written as sums of
Fej\'{e}r graph integrals, by computing bounds of the form
\begin{equation*}
T^{d\alpha _{k}(z)}
\end{equation*}%
for Fej\'{e}r graph integrals intervening in the cumulants of order $k$.
Since the typical CLT normalization is $T^{d/2}$, it will suffice then
establishing \textbf{\textquotedblleft cumulant inequalities"}
\begin{equation*}
d\alpha _{k}(z)<kd/2\Leftrightarrow \boxed{ \alpha _k( z) < k /2}
\end{equation*}%
where $\alpha _{k}(z)$ is the exponent appearing in the expansion of the $k$%
-th cumulant. In fact, this may be strengthened (cf. Theorem \ref{Sz}) to
include $z$ satisfying the equality $\alpha _{k}(z)=k/2$, if $k\geq 3$.

In conclusion, establishing normality is reduced to computing the functions $%
\alpha _k( z), k \geq 3,$ i.e., to solving a sequence of {graph breaking
problems}.

\begin{xmpl}
\label{exagraph} The general structure of the intervening graphs for the $k$%
'th cumulant of sums $S_T$ of the $m$'th Appell polynomial of a Gaussian
sequence is provided by graphs belonging to the set $\Gamma(m,k)$ of all
connected graphs with no loops over $k$ vertices, each of degree $m$ (see
\cite{a}). Let $z=p^{-1}$ denote the integrability exponent of the spectral
density.

The cumulant inequality corresponding to the \textbf{\textquotedblleft
maximal breaking=MB"} of the $k$'th cumulant graph, which typically yields a
facet of the power counting polytope (PCP) (at least for $k=2$) is given in
this example by:

\begin{equation*}
\alpha _k( z)= \sum_e (z_e-1) + V =\frac{k m }{2} (z-1) + k \leq \frac{k }{2}
\Leftrightarrow \frac{1}{ m } \leq 1-z
\end{equation*}

At the limiting point $1-z=\frac{1}{ m }$, the cumulant exponents are $%
\alpha _k(z)= \frac{ k}{ 2 }$, ensuring negligibility for $k \geq 3$.

The discrete and continuous case may be unified here by asking for
integrability at $z=1-\frac{1}{ m }$, since in the discrete case the
extension to smaller values of $z$ is trivial.
\end{xmpl}


\begin{xmpl}
\label{exaquadr} A similar analysis holds in the case of cumulants of
quadratic forms in Appell polynomials $P_{m,n}(X_t,X_s)$. Note that while
the number of graphs intervening increases considerably, the number of
extremal points of the PCP is just $4$ -- see Figure \ref{fpol}.

\begin{figure}[h!]
\begin{center}
\leavevmode
\resizebox{.82 \textwidth}{!}
{\ \includegraphics{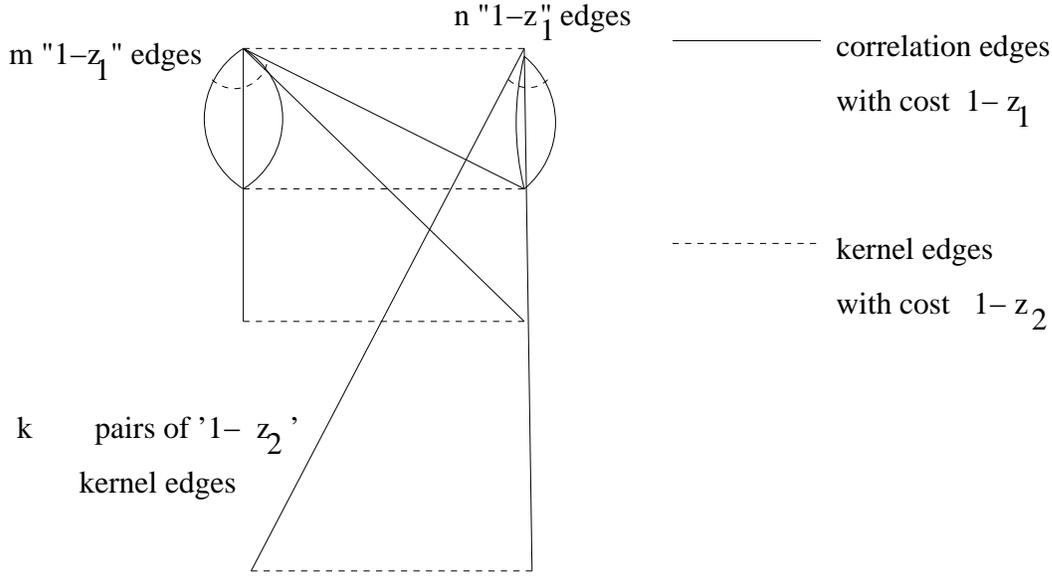}}
\end{center}
\caption{\textbf{The graphs appearing in the expansion of cumulants of
quadratic forms. Here k=4, m=5, n=4. The figure displays only some of the
k(m+n)/2=18 correlation edges.}}
\label{fgr}
\end{figure}

The graphs intervening are -- see Figure \ref{fgr} -- all the graphs
belonging to the set $\Gamma(m,n,k)$ of all connected bipartite graphs with
no loops whose vertex set consists of $k$ pairs of vertices. The ``left"
vertex of each pair arises out of the first $m$ terms $:X_{t_1},
...,X_{t_m}: $ in the diagram formula, and the ``right" vertex of each pair
arises out of the last $n$ terms $:X_{s_1}, ...,X_{s_n}:$ The edge set
consists of:

\begin{enumerate}
\item $k$ ``kernel edges" pairing each left vertex with a right vertex. The
kernel edges will contribute below terms involving the function $b(\lambda)$.

\item A set of ``correlation edges", always connecting vertices in different
rows, and contributing below terms involving the function $f(\lambda)$).
They are arranged such that each left vertex connects to $m$ and each right
vertex connects to $n$ such edges, yielding a total of $k(m+n)/2$
correlation edges.
\end{enumerate}

Thus, the $k$ ``left vertices" are of degree $m+1$, and the other $k$
vertices are of degree $n+1$. (The ``costs $1-z_1, 1-z_2$'' mentioned in
Figure \ref{fgr} refer to (\ref{borne})).

The PCP domain in the discrete case (which is precisely the convergence
domain of the integrals defining the limiting variance), is indicated below,
when $m < n$, in terms of the integrability indices $z=(z_1, z_2)$ of $f$
and $b$ (i.e. $f \in L_{p_1}, b \in L_{p_2}$).

\begin{figure}[h!]
\begin{center}
\leavevmode
\resizebox{.82 \textwidth}{!}
{\ \includegraphics{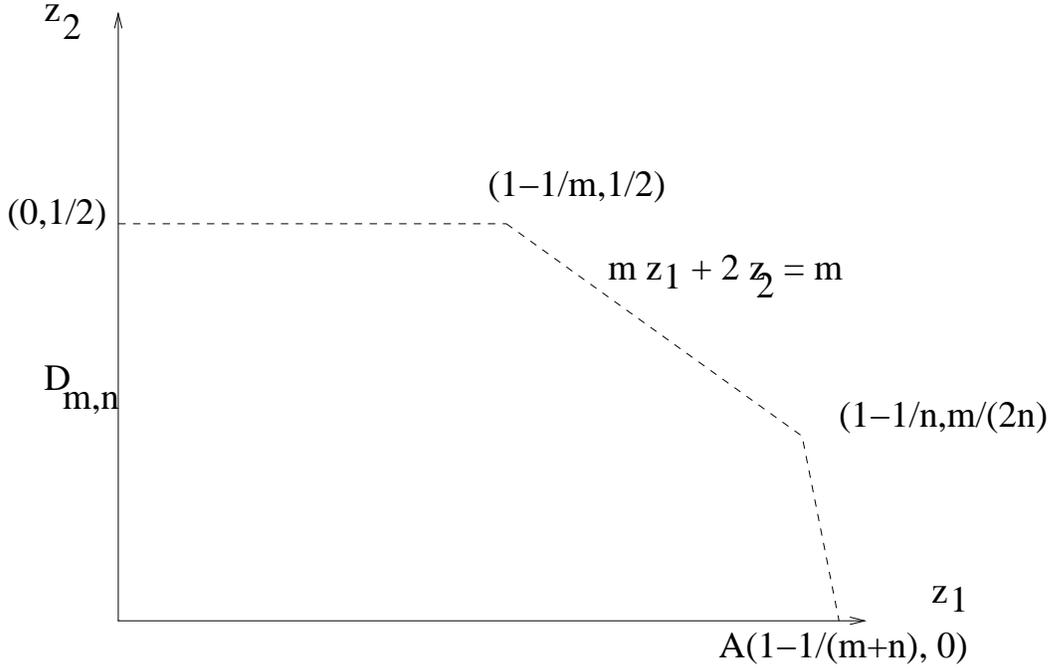}}
\end{center}
\caption{\textbf{The domain of the central limit theorem, discrete case}}
\label{fpol}
\end{figure}

When $m<n$, there are only three segments on the undominated boundary of the
PCP, connecting respectively the extremal points $(A,B)$, $(B,C)$ and $%
(C,D)( $with coordinates $A(1-1/(m+n),0),$ $B(1-1/n,m/(2n)),$ $C(1-1/m,1/2),$
$D(0,1/2))$, and correspond respectively to the breakings indicated below:
\begin{equation*}
\begin{cases}
\text{total breaking } & 2k-\frac{k(m+n)}{2}(1-z_{1})-k(1-z_{2})\leq \frac{k%
}{2} \\
& \Leftrightarrow \frac{3}{2}\leq \frac{(m+n)}{2}(1-z_{1})+(1-z_{2}) \\
\text{breaking all $z_{2}$ and the left $z_{1}$ edges} & \max_{k}\frac{k+1}{k%
}-m(1-z_{1})-(1-z_{2})\leq \frac{1}{2} \\
& \Leftrightarrow 1\leq m(1-z_{1})+(1-z_{2}) \\
\text{breaking all $z_{2}$ edges } & \max_{k}\frac{2}{k}-(1-z_{2})=\leq
\frac{1}{2} \\
& \Leftrightarrow \frac{1}{2}\leq z_{2}%
\end{cases}%
\end{equation*}

In the continuous case, the domain is just the lower segment between the
points $A$ and $B$ in figure \ref{fpol}.
\end{xmpl}

\subsection{Proof of Theorem \protect\ref{Sz}, parts a), c)}

We turn now to theorem \ref{Sz} a), c), generalizing Theorem 2 and Corollary
1 of \cite{a}.

For part a), let us apply the H\"older-Young-Brascamp-Lieb inequality with
optimally chosen integrability parameters $s_v^{-1}$:
\begin{eqnarray}
|J_T(M, f_1,...,f_E)|\leq K \prod_{v=1}^V |\!|\Delta_T|\!|_{s_v^{-1}}
\prod_{e=1}^E |\!|f_e|\!|_{z_e^{-1}} \\
\leq K^{\prime}T^{d \, \sum_{v=1}^V (1- s_v^{-1})} \prod_{e=1}^E
|\!|f_e|\!|_{z_e^{-1}}
\end{eqnarray}
under the constraint that $(s_1, ..., s_V, z_1, ..., z_E)$ satisfy the power
counting conditions, and where we used the kernel estimate
\begin{equation*}
|\!|\Delta _{T}(\lambda )|\!|_{s^{-1}} \leq C_{s} T^{d(1-s)}, \forall s \in
[0,1)
\end{equation*}
(see Appendix B).

The optimization problem for $s_v$ in the discrete case:
\begin{eqnarray*}
\min_{s_1, ..., s_V} d \; \sum_{v=1}^V (1- s_v^{-1}), \text{ where } (s_1,
..., s_V, z_1, ..., z_E) \in \text{ PCP},
\end{eqnarray*}
has the same constraints as Lemma 2 in \cite{a}, except that the objective
is multiplied by $d$. Hence, in the torus case, the exponent is simply $d$
times the one dimensional exponent of Theorem 1, \cite{a}.

In the continuous case $S=\mathbb{R}^{d}$, we note first that when $%
\sum_{e}z_{e}\leq C$ the result follows, just as Theorem \ref{Sz} b), from
the H\"{o}lder-Young-Brascamp-Lieb \ inequality, while in the other case $%
\sum_{e}z_{e}\geq C$, the extra constraint $\sum_{v}s_{v}=E-\sum_{e}z_{e}$
yields the one-dimensional exponent as $\alpha _{M}(z)=co(M)+\sum_{e}z_{e}-C$%
. 

For part c), we approximate our functions by   continuous, bounded
functions  of bounded support,  for which conditions (3.3) of part
b) hold. For these approximants, it follows from the
convergence $\frac{J_{T}(M)}{T^{d\;co(M)}}\rightarrow k_{M}\;\mathcal{J}%
(M^{\ast },f_{1},...,f_{E})$  that $J_{T}(M)=o(T^{d\;co(M)+a)})$, $%
\forall a>0.$  The result follows then for the functions from
${{\mathbf{L}}}_{p}$ spaces, by the definition of these spaces. Note
that in the discrete case, the same argument was applied based on
trigonometric polynomials -- see \cite{a}, proof of Corollary 1.





\section{Applications \label{s:app}}

\subsection{Central limit theorems for bilinear forms of moving averages
\label{ss:lin}}

We assume below that our stationary random field $X_{t},\ t\in I$ admits a
representation as a linear/moving averages random field.

For discrete parameter it means that
\begin{equation}
X_{t}=\sum_{u\in \mathbb{Z}^{d}}\hat{a}(t-u)\xi ({u}),\;\sum_{u\in \mathbb{Z}%
^{d}}\hat{a}^{2}(u)<\,\infty ,~t\in \mathbb{Z}^{d},  \label{nn1}
\end{equation}%
where $\xi (u)$, $u\in \mathbb{Z}^{d}$, are independent random variables
indexed by $\mathbb{Z}^{d}$ with ${\mathbb{E}}\xi (0)=0$ and such that $%
E\left\vert \xi (0)\right\vert ^{k}\leq c_{k}<\infty ,~k=1,2,\ldots .$ In
this case
\begin{equation}
c_{k}(t_{1},...,t_{k})=cum_{k}\{X_{t_{1}},...,X_{t_{k}}\}=d_{k}\sum%
\limits_{s\in \mathbb{Z}^{d}}\prod_{j=1}^{k}\hat{a}(t_{j}-s),  \label{nn2}
\end{equation}%
where $d_{k}$ is the $k$'th cumulant of $\xi (0).$

For continuous parameter we assume that
\begin{equation}
X_{t}=\int_{u\in {\mathbb{R}^{d}}}\hat{a}(t-u)\xi ({du}),\quad t\in {\mathbb{%
R}^{d}},  \label{nn3}
\end{equation}
with a square-integrable kernel $\hat{a}(t)$, $t\in {I}$, with respect to a
independently scattered random measure with finite second moment, that is a
homogeneous random measure $\xi (A),\ A\subset {\mathbb{R}^{d}}$, with
finite second moments and independent values over disjoint sets (see, for
instance, Rajput and Rosinski \cite{RajputRosinski1989} or Kwapien and
Woyczynski \cite{KwapienWoyczynski1992}). That is, for each Borel $A,~\xi
(A) $ is an infinitely divisible random variable whose cumulant function can
be written as
\begin{equation}
\kappa (z)=\log Ee^{iz\xi (A)}=izm_{0}(A)-\frac{1}{2}z^{2}m_{1}(A)+\int_{%
\mathbb{R}}\left( e^{izx}-1-iz\tau (x)\right) Q(A,dx),  \label{nn4}
\end{equation}
where $m_{0}$ is a signed measure, $m_{1}$ is a positive measure, $Q(A,dx)$
(for fixed $A$) is a measure on ${\mathbb{R}^{1}}$ without atoms at $0$,
such that $\int_{\mathbb{R}}\min \left\{ 1,\left| x\right| ^{2}\right\}
Q(A,dx)<\infty $, and where $\tau (x)=x$ if $\left| x\right| \leq 1,$ and $%
\tau (x)=x/\left| x\right| ,$ if $\left| x\right| >1.$

For example, if $I=\mathbb{R}$, then $\xi (A)$ is a set indexed L\'{e}vy
process with finite second moments and stationary intensity proportional to
the Lebesgue measure.

We also assume that $Q$ factorizes as $Q(A,dx)=M(A)W(dx),$ where $M(A)$ is a
$\sigma -$finite measure, and $W$ is some L\'{e}vy measure on $\mathbb{R}%
^{1},$ such that for some $\varepsilon >0$ and $\lambda >0$%
\begin{equation*}
\int_{(-\varepsilon ,\varepsilon )}e^{\lambda u}W(du)<\infty .
\end{equation*}

This implies that
\begin{equation*}
\int_{\mathbb{R}}\left\vert u\right\vert ^{k}W(du)<\infty ,~k\geq 2,
\end{equation*}
and that the cumulant function $\kappa (z)$ is analytical in a neighborhood
of $0.$

Necessary and sufficient conditions of existence of the integral (as limit
in probability of integrals of simple functions)
\begin{equation*}
\int_{A}f(s)d\xi (s)
\end{equation*}
can be found in \cite{RajputRosinski1989}. Note that, for $d=1,$ an
integrals become integrals with respect to L\'{e}vy process $L(t),~t\in
\mathbb{R}^{1}$, and $\kappa \left( z\right) =\log Ee^{izL(1)}.$

For Lebesgue measures $m_{0},~m_{1}$ and $Q,$ one can prove (by using
product integration) that
\begin{equation}
\log E\exp \left\{ i(z_{1}X_{t_{1}}+\cdots +z_{k}X_{t_{k}}\right\}
=\int\limits_{\mathbb{R}^{d}}\kappa \left( \sum_{j=1}^{k}z_{j}\hat{a}%
(t_{j}-s)\right) ds  \label{nn5}
\end{equation}
if $\hat{a}\in L_{1}\cap L_{2}.$ From (\ref{nn5}) if can be seen that random
field (\ref{nn3}) is homogeneous in a strict sense.

We assume that $m_{1}=0,$ that is, ${\mathbb{E}}\xi (I_{1})=0,$ then the
last formula holds for $\hat{a}\in L_{2}.$

We obtain that
\begin{equation}
c_{k}(t_{1},...,t_{k})=cum_{k}\{X_{t_{1}},...,X_{t_{k}}\}=d_{k}\int_{\mathbb{%
R}^{d}}\prod_{j=1}^{k}\hat{a}(t_{j}-s)ds,  \label{nn6}
\end{equation}
where $d_{k}$ is the k'th cumulant of $\xi (I_{1})$ with $I_{1}$ being the
unit rectangle, that is $d_{k}=\kappa ^{k}(0)/i^{k},~k\geq 2.$

We assume from now on that ${\mathbb{E}}\xi (I_{1})=0$, and use the same
notation for both discrete and continuous cases

\begin{equation}
X_{t}=\int_{u\in I}\hat{a}(t-u)\xi ({du}),\quad t\in I,  \label{nn7}
\end{equation}
where $I=\mathbb{Z}^{d}$ in the discrete case and $I=\mathbb{R}^{d}$ in the
continuous case.

For various conditions which ensure that (\ref{nn7}) is well-defined, see,
for example, Anh, Heyde and Leonenko \cite{AHL}, p. 733, and references
therein.

\textbf{Note:} By choosing an appropriate ``Green function'' $\hat{a}(t)$,
this very general class of processes includes the solutions of many
interesting differential equations with random noise $\xi ({du})$, like, for
example, generalized Ornstein-Uhlenbeck processes in $\mathbb{R}$ \cite{AHL}.

We will assume here that all moments for our stationary field $X_{t}$ exist.

The advantage of the linear representation assumption (\ref{nn7}) and (\ref%
{nn4}), (\ref{nn5}) is the explicit representation of cumulants -- see for
example Theorem 2.1 of \cite{AHL}:
\begin{equation}
c_{k}(t_{1},...,t_{k})=d_{k}\int_{s\in I}\prod_{j=1}^{k}\hat{a}%
(t_{j}-s)\;\nu (ds),  \label{nn8}
\end{equation}
where $d_{k}$ is the k'th cumulant of $\xi (I_{1})$ with $I_{1}$ being the
unit rectangle.

In the spectral domain, we get
\begin{equation}
f_{k}(\lambda _{1},...,\lambda _{k-1})=d_{k}\;a(-\sum_{i=1}^{k-1}\lambda
_{i})\;\prod_{i=1}^{k-1}a(\lambda _{i})=\prod_{i=1}^{k}a(\lambda _{i})\delta
(\sum_{j=1}^{k}\lambda _{j}).  \label{nn9}
\end{equation}

For $k=2$, we will denote the spectral density by $f(\lambda )=f_{2}(\lambda
)=d_{2}a(\lambda )a(-\lambda )$.

We can formulate now a central limit theorems for quadratic functional of a
linear field, which is a generalization of the results of Giraitis and
Surgailis \cite{giraitis:surgailis:1986} and Giraitis and Taqqu \cite%
{giraitis:taqqu:1997} (see also references therein) for discrete time
processes. This next theorem follows from the results of Sections 3 and 4
(the proof is almost identical to the proof of Theorem 4 of Avram \cite{a},
the expression for the variance can be obtained by direct computations).

\begin{thrm}
\label{linbilin} Let $X_{t}=\int_{u\in {\ I}}\hat{a}(t-u)\xi ({\ du}),\quad
t\in {\ I,}$ be a linear random field with a square integrable kernel $\hat{a%
}(t),t\in I,$ and a random measure $\xi (du)$ admitting all moments and let%
\begin{equation*}
Q_{T}=Q_{T}^{(1,1)}=\int_{t,s\in I_{T}}\left[ X_{t}X_{s}-{\mathbb{E}}%
X_{t}\,X_{s}\right] \hat{b}(t-s)\nu (ds)\nu (dt).
\end{equation*}%
We assume that $f(\lambda )=(2\pi )^{d}\left\vert a(\lambda )\right\vert
^{2}\in \mathbf{L}_{p}$ and $b(\lambda )\in \mathbf{L}_{q},$ and in the
continuous case we assume also that $b(\lambda )\in \mathbf{L}_{q}\cap
L_{1}. $

Assume that:
\begin{equation*}
\frac{1}{p}+\frac{1}{q}\leq \frac{1}{2}.
\end{equation*}%
Then, the central limit theorem holds:
\begin{equation*}
\lim T^{-d/2}Q_{T}\rightarrow N(0,\sigma ^{2}),\text{ \ \ }T\rightarrow
\infty ,
\end{equation*}%
where
\begin{equation*}
\sigma ^{2}:=2(2\pi )^{d}d_{2}^{2}\int_{S}b^{2}(\lambda )f^{2}(\lambda
)d\lambda +(2\pi )^{d}d_{4}\left( \int_{S}b(\lambda )f(\lambda )d\lambda
\right) ^{2}
\end{equation*}%
where $d_{k}$ is the k'th cumulant of $\xi (I_{1})$ with $I_{1}$ being the
unit rectangle, that is

$d_{2}={\mathbb{E}}\xi (I_{1})^{2},d_{4}={\mathbb{E}}(\xi (I_{1})^{4})-2[{%
\mathbb{E(}}\xi (I_{1})^{2})]^{2}$ in the continuous case and $d_{k}$ is the
k'th cumulant of $\xi (0)$ in the discrete case, that is $d_{2}={\mathbb{E}}%
\xi (0)^{2},$ $d_{4}={\mathbb{E(}}\xi (0)^{4})-2[{\mathbb{E(}}\xi
(0)^{2})]^{2}.$
\end{thrm}

\subsection{Minimum contrast estimation based on the Whittle contrast
function}

The class of Whittle estimators is the most popular in applications (see
Whittle \cite{Wh51}, \cite{Wh53}, Giraitis and Surgailis \cite%
{GiraitisSurgailis1990}, Fox and Taqqu \cite{ft}, Heyde and Gay \cite{hg},
\cite{hg1}, Heyde \cite{Heyde}, Gao, Anh and Heyde \cite{GAH2002}, Leonenko
and Sakhno \cite{LS2006}, see also the references therein).

In what follows we will consider continuous time linear processes $(d=1)$
whose spectral densities of all orders exist and admit the representation of
the form (\ref{nn9}).

We begin with the following assumption.

\textbf{A.I}. Let $X_{t},$ $t\in I_{T}=\left[ -\frac{T}{2},\frac{T}{2}\right]
,$ be an observation of a real-valued measurable stationary linear process $%
X_{t},$ $t\in \mathbb{R}^{1}$, with zero mean and the family of spectral
densities (\ref{nn9}). Let $a\left( \lambda \right) =a\left( \lambda ;\theta
^{(1)}\right) ,d_{k}=d_{k}\left( \theta ^{(2)}\right) ,$ that is, $%
f_{2}(\lambda )=f\left( \lambda ,\theta \right) ,$ $\lambda \in \mathbb{R}%
^{1},$ $\theta =\left( \theta ^{(1)},\theta ^{(2)}\right) ,$ $\theta \in
\Theta \subset \mathbb{R}^{m},$ where $\Theta $ is a compact set, and the
true value of the parameter $\theta _{0}\in int\Theta ,$ the interior of $%
\Theta .$ Suppose further that $f\left( \lambda ;\theta _{1}\right)
\not\equiv f\left( \lambda ;\theta _{2}\right) $ for $\theta _{1}\neq \theta
_{2},$ almost everywhere in $\mathbb{R}^{1}$ with respect to the Lebesgue
measure.

Consider the Whittle contrast process (or objective function)
\begin{equation}
U_{T}\left( \theta \right) =\frac{1}{4\pi }\int_{\mathbb{R}^{1}}\left( \log
f\left( \lambda ;\theta \right) +\frac{I_{T}\left( \lambda \right) }{f\left(
\lambda ;\theta \right) }\right) w\left( \lambda \right) d\lambda ,
\label{(2W)}
\end{equation}%
where $I_{T}\left( \lambda \right) $ is the periodogram of the second order
\begin{equation}
I_{T}\left( \lambda \right) =\frac{1}{2\pi T}\left\vert
\int_{I_{T}}X_{t}e^{-it\lambda }dt\right\vert ^{2},\text{ \ }\lambda \in
\mathbb{R}^{1},
\end{equation}%
and $w\left( \lambda \right) $ is a symmetric about $\lambda =0$ function
such that all considered integrals are well defined and which will satisfy
some conditions given below; in some cases we can choose $w\left( \lambda
\right) =\frac{1}{1+\lambda ^{2}}.$

Introduce the Whittle contrast function
\begin{equation}
K\left( \theta _{0};\theta \right) =\frac{1}{4\pi }\int_{\mathbb{R}%
^{1}}\left( \frac{f\left( \lambda ;\theta _{0}\right) }{f\left( \lambda
;\theta \right) }-1-\log \frac{f\left( \lambda ;\theta _{0}\right) }{f\left(
\lambda ;\theta \right) }\right) w\left( \lambda \right) d\lambda .
\label{(4)}
\end{equation}

To state the result on consistency of the minimum contrast estimator based
on the contrast process (\ref{(2W)}) we will need the following conditions
on the spectral density $f\left( \lambda ;\theta \right) $ and the weight
function $w\left( \lambda \right) .$

\noindent\textbf{A}.\textbf{II. }$f\left( \lambda ;\theta _{0}\right)
w\left( \lambda \right) \frac{1}{\mathbf{\ }f\left( \lambda ;\theta \right) }%
\in L_{1}\left( \mathbb{R}^{1}\right) \cap L_{2}\left( \mathbb{R}^{1}\right)
,\quad \forall \theta \in \Theta .$\textbf{\ }

\noindent\textbf{A.III.} There exists a function $v\left( \lambda \right)
,\lambda \in \mathbb{R}^{1},$ such that

(i) the function $h\left( \lambda ;\theta \right) =v\left( \lambda \right)
\frac{1}{\mathbf{\ }f\left( \lambda ;\theta \right) }$ is uniformly
continuous in $\mathbb{R}^{1}\times \Theta $;

(ii) $f\left( \lambda ;\theta _{0}\right) \frac{w\left( \lambda \right) }{%
v\left( \lambda \right) }\in L_{1}\left( \mathbb{R}^{1}\right) \cap
L_{2}\left( \mathbb{R}^{1}\right) .$

\begin{thrm}
\label{Wh1} Let the assumptions A.I to A.III be satisfied. Then the function
$K\left( \theta _{0};\theta \right) $ defined by (\ref{(4)}) is the contrast
function for the contrast process $U_{T}\left( \theta \right) $ defined by (%
\ref{(2W)}). The minimum contrast estimator $\widehat{\theta }_{T}$ defined
as
\begin{equation}
\widehat{\theta }_{T}=\underset{\theta \in \Theta }{\arg \min }U_{T}\left(
\theta \right)  \label{(5)}
\end{equation}%
is a consistent estimator of the parameter $\theta ,$ that is, $\widehat{%
\theta }_{T}\rightarrow \theta _{0}$ in $P_{0}$-probability as $T\rightarrow
\infty .$
\end{thrm}

The above theorem can be obtained as a consequence of a more general result
by Leonenko and Sakhno \cite{LS2006} (Theorem 3), one needs just to rewrite\
for the case of linear processes the corresponding conditions on spectral
densities, which become of much simpler form.

Next set of assumptions (in addition to the above ones) is needed to state
the result on asymptotic normality of the estimator (\ref{(5)}).

\noindent\textbf{A.IV. \ }The function $\frac{1}{\mathbf{\ }f\left( \lambda
;\theta \right) }$ is twice differentiable in a neighborhood of the point $%
\theta _{0}$ and

\noindent(i) $f\left( \lambda ;\theta _{0}\right) w\left( \lambda \right)
\frac{\partial ^{2}}{\partial \theta _{i}\partial \theta _{j}}\frac{1}{%
\mathbf{\ }f\left( \lambda ;\theta \right) }$ $\in $ $L_{1}\left( \mathbb{R}%
^{1}\right) \cap L_{2}\left( \mathbb{R}^{1}\right) ,$ $i,j=1,...,m,$ $\theta
\in \Theta ; $

\noindent (ii) $f\left( \lambda ;\theta _{0}\right) \in \mathbf{L}_{p}\left(
\mathbb{R}^{1}\right) ,$ $w\left( \lambda \right) \frac{\partial }{\partial
\theta _{i}}\frac{1}{\mathbf{\ }f\left( \lambda ;\theta \right) }$ $\in $ $%
\mathbf{L}_{q}\left( \mathbb{R}^{1}\right) ,$\

for some $p,$ $q$ such that $\frac{1}{p}+\frac{1}{q}\leq \frac{1}{2}$, $%
i=1,...,m,$ $\theta \in \Theta ;$

\noindent(iii) $T^{1/2}\int_{\mathbb{R}^{1}}(EI_{T}\left( \lambda \right)
-f\left( \lambda ;\theta _{0}\right) )w\left( \lambda \right) \frac{\partial
}{\partial \theta _{i}}\frac{1}{\mathbf{\ }f\left( \lambda ;\theta \right) }%
\ d\lambda \rightarrow 0$ \ as $T\rightarrow \infty ,$

for all $i=1,...,m,\theta \in \Theta ;$

\noindent (iv) the second order derivatives $\frac{\partial ^{2}}{\partial
\theta _{i}\partial \theta _{j}}\frac{1}{\mathbf{\ }f\left( \lambda ;\theta
\right) },i=1,...,m,$ are continuous in $\theta .$

\noindent \textbf{A.V. }The matrices $W_{1}\left( \theta \right) =\left(
w_{ij}^{(1)}\left( \theta \right) \right) _{i,j=1,...,m}$, $W_{2}\left(
\theta \right) =\left( w_{ij}^{(2)}\left( \theta \right) \right)
_{i,j=1,...,m}$, $V\left( \theta \right) =\left( v_{ij}\left( \theta \right)
\right) _{i,j=1,...,m}$ are positive definite, where
\begin{equation}
w_{ij}^{(1)}\left( \theta \right) =\frac{1}{4\pi }\int_{\mathbb{R}%
^{1}}w\left( \lambda \right) \frac{\partial }{\partial \theta _{i}}\log
f\left( \lambda ;\theta \right) \frac{\partial }{\partial \theta _{j}}\log
f\left( \lambda ;\theta \right) d\lambda ,
\end{equation}
\begin{equation}
w_{ij}^{(2)}\left( \theta \right) =\frac{1}{4\pi }\int_{\mathbb{R}%
^{1}}w^{2}\left( \lambda \right) \frac{\partial }{\partial \theta _{i}}\log
f\left( \lambda ;\theta \right) \frac{\partial }{\partial \theta _{j}}\log
f\left( \lambda ;\theta \right) d\lambda .
\end{equation}

\begin{equation}
v_{ij}\left( \theta \right) =\frac{1}{8\pi }\frac{d_{4}}{d_{2}^{2}}\int_{%
\mathbb{R}^{1}}w\left( \lambda \right) \frac{\partial }{\partial \theta _{i}}%
\log f\left( \lambda ;\theta \right) d\lambda \int_{\mathbb{R}^{1}}w\left(
\lambda \right) \frac{\partial }{\partial \theta _{j}}\log f\left( \lambda
;\theta \right) d\lambda .\text{ }
\end{equation}

\begin{thrm}
\label{Wh2} Let the assumptions A.I to A.V be satisfied. Then as $%
T\rightarrow \infty $%
\begin{equation*}
T^{1/2}\left( \widehat{\theta _{T}}-\theta _{0}\right) \overset{\mathcal{D}}{%
\rightarrow }\mathcal{N}_{m}\left( 0,W_{1}^{-1}\left( \theta _{0}\right)
(W_{2}\left( \theta _{0}\right) +V\left( \theta _{0}\right)
)W_{1}^{-1}\left( \theta _{0}\right) \right) ,
\end{equation*}
where $N_{m}\left( \cdot ,\cdot \right) $ denotes the $m$-dimensional
Gaussian law.
\end{thrm}

Reasonings for the proof of Theorems \ref{Wh1}, \ref{Wh2} are given in the
next section.

Comparing the above theorem with a more general result stated in \cite%
{LS2006}, one can see that the set of conditions for the case of linear
processes becomes of much simpler form, but the most important improvement
is in condition A.IV(ii), which was achieved due to the application of the
Theorem \ref{linbilin} (see the proof). Note that corresponding condition
for the case of general processes, formulated in \cite{LS2006},
unfortunately, is difficult to check in general situation.

\begin{rmrk}
Condition A.IV(iii) will hold, e.g., if $f\left( \lambda ;\theta \right) $
is differentiable with respect to $\lambda $ and
\begin{equation*}
\int_{\mathbb{R}^{1}}f_{\lambda }^{\prime }\left( \lambda ;\theta
_{0}\right) )w\left( \lambda \right) \frac{\partial }{\partial \theta _{i}}%
\frac{1}{\mathbf{\ }f\left( \lambda ;\theta \right) }d\lambda <\infty ,
\end{equation*}
or under any conditions which assure
\begin{equation*}
\int_{\mathbb{R}^{1}}|f\left( \lambda +h;\theta _{0}\right) -f\left( \lambda
;\theta _{0}\right) |w\left( \lambda \right) \frac{\partial }{\partial
\theta _{i}}\frac{1}{\mathbf{\ }f\left( \lambda ;\theta \right) }d\lambda
\leq C|h|^{a},
\end{equation*}
for $a>\frac{1}{2}$ and $C$ being a constant.
\end{rmrk}

\textbf{Example.} Estimation of fractional Riesz-Bessel motion (FRBM) (see
Appendix A for details and definition of FRBM in non-Gaussian case). Let $%
X_{t},$ $t\in \mathbb{R}^{1},$ be a non-Gaussian Riesz-Bessel stationary
motion, that is a stationary linear process with the spectral density of the
form
\begin{equation}
f\left( \lambda \right) =f\left( \lambda ,\theta \right) =\frac{c}{%
\left\vert \lambda \right\vert ^{2\alpha }\left( 1+\lambda ^{2}\right)
^{\gamma }},\quad \lambda \in \mathbb{R}^{1},
\end{equation}%
where the unknown vector parameter $\theta =\left( \gamma ,\alpha ,c\right)
^{\prime }\in \Theta ,\,\Theta $ being a compact subset of $\left[ \frac{1}{2%
},\infty \right) \times \left( 0,\frac{1}{2}\right) \times \left( 0,\infty
\right) .$ Note that the index $\alpha $ determines the long-range
dependence of FRBM, and the parameter $\gamma $ is another fractal index
connected to Hausdorff dimension of paths of the stochastic process. Note
that procedure of discretazation leads to the loss of information of one
parameter $\gamma $, which is important for applications in both turbulence
and finance theory. That is why a direct method of estimation of both
parameters from continuous data looks appropriate.

For this model we can choose the weight function $w\left( \lambda \right) =%
\frac{1}{1+\lambda ^{2}},$ $\lambda \in \mathbb{R}^{1},$ to satisfy the
conditions needed for consistency of the estimator (\ref{(5)}), that is, for
Theorem \ref{Wh1} to hold. However, to satisfy all the conditions needed for
Theorem \ref{Wh2} we choose the weight function $w\left( \lambda \right) =%
\frac{\lambda ^{2b}}{\left( 1+\lambda ^{2}\right) ^{a}},$ $\lambda \in
\mathbb{R}^{1},$ where $a$ and $b$ satisfy the restrictions:$\{b>1\}\wedge
\{a>b+2\}\wedge \{a>A+2\},$where we have denoted by $A$ the length of the
finite interval carrying the admissible values of the parameter $\gamma $.
With such a choice of the weight function we have the convergence
\begin{equation*}
T^{1/2}\left( \widehat{\theta _{T}}-\theta _{0}\right) \overset{\mathcal{D}}{%
\rightarrow }N_{3}\left( 0,W_{1}^{-1}\left( \theta _{0}\right) (W_{2}\left(
\theta _{0}\right) +V\left( \theta _{0}\right) )W_{1}^{-1}\left( \theta
_{0}\right) \right) \text{ \ as \ }T\rightarrow \infty ,
\end{equation*}%
where the elements of the matrices $W_{1}$ and $W_{2}$\ are of the following
form:
\begin{eqnarray*}
w_{{}11}^{(1\vee 2)} &=&\frac{1}{4\pi }\int_{\mathbb{R}^{1}}w^{1\vee
2}\left( \lambda \right) \left( \ln \left( 1+\lambda ^{2}\right) \right)
^{2}d\lambda ; \\
w_{{}22}^{(1\vee 2)} &=&\frac{1}{4\pi }\int_{\mathbb{R}^{1}}w^{1\vee
2}\left( \lambda \right) \left( \ln \left( \lambda ^{2}\right) \right)
^{2}d\lambda ; \\
w_{{}33}^{(1\vee 2)} &=&\frac{1}{4\pi }c_{0}^{-2}\int_{\mathbb{R}%
^{1}}w^{1\vee 2}\left( \lambda \right) d\lambda ;
\end{eqnarray*}%
\begin{eqnarray*}
w_{{}12}^{(1\vee 2)} &=&w_{{}21}^{(1\vee 2)}=\frac{1}{4\pi }\int_{\mathbb{R}%
^{1}}w^{1\vee 2}\left( \lambda \right) \ln \left( 1+\lambda ^{2}\right) \ln
\left( \lambda ^{2}\right) d\lambda ; \\
w_{{}13}^{(1\vee 2)} &=&w_{{}31}^{(1\vee 2)}=-\frac{1}{4\pi }c_{0}^{-1}\int_{%
\mathbb{R}^{1}}w^{1\vee 2}\left( \lambda \right) \ln \left( 1+\lambda
^{2}\right) d\lambda ; \\
w_{{}23}^{(1\vee 2)} &=&w_{{}32}^{(1\vee 2)}=-\frac{1}{4\pi }c_{0}^{-1}\int_{%
\mathbb{R}^{1}}w^{1\vee 2}\left( \lambda \right) \ln \left( \lambda
^{2}\right) d\lambda ;
\end{eqnarray*}%
\begin{eqnarray*}
v_{11} &=&\frac{1}{8\pi }\frac{d_{4}}{d_{2}^{2}}\left( \int_{\mathbb{R}%
^{1}}w\left( \lambda \right) \ln \left( 1+\lambda ^{2}\right) d\lambda
\right) ^{2}; \\
v_{22} &=&\frac{1}{8\pi }\frac{d_{4}}{d_{2}^{2}}\left( \int_{\mathbb{R}%
^{1}}w\left( \lambda \right) \ln \left( \lambda ^{2}\right) d\lambda \right)
^{2}; \\
v_{33} &=&\frac{1}{8\pi }\frac{d_{4}}{d_{2}^{2}}c_{0}^{-2}\left( \int_{%
\mathbb{R}^{1}}w\left( \lambda \right) d\lambda \right) ^{2};
\end{eqnarray*}%
\begin{eqnarray*}
v_{12} &=&v_{21}=\frac{1}{8\pi }\frac{d_{4}}{d_{2}^{2}}\int_{\mathbb{R}%
^{1}}w\left( \lambda \right) \ln \left( 1+\lambda ^{2}\right) d\lambda \int_{%
\mathbb{R}^{1}}w\left( \lambda \right) \ln \left( \lambda ^{2}\right)
d\lambda ; \\
v_{13} &=&v_{31}=-\frac{1}{8\pi }\frac{d_{4}}{d_{2}^{2}}c_{0}^{-1}\int_{%
\mathbb{R}^{1}}w\left( \lambda \right) \ln \left( 1+\lambda ^{2}\right)
d\lambda \int_{\mathbb{R}^{1}}w\left( \lambda \right) d\lambda ; \\
v_{23} &=&v_{32}=-\frac{1}{8\pi }\frac{d_{4}}{d_{2}^{2}}c_{0}^{-1}\int_{%
\mathbb{R}^{1}}w\left( \lambda \right) \ln \left( \lambda ^{2}\right)
d\lambda \int_{\mathbb{R}^{1}}w\left( \lambda \right) d\lambda .
\end{eqnarray*}%
In the above formulae we mean that the weight function $w\left( \lambda
\right) $ is involved to the expressions for $w_{{}ij}^{(1)}$ in the 1st
power and to the expressions for $w_{{}ij}^{(2)}$ in the 2d power. From the
above formulae we see that the covariance matrix of the limiting normal law
has the charming feature that it appears not depending on the values $\alpha
_{0}$ and $\gamma _{0}$.

\begin{rmrk}
Continuous version of Gauss-Whittle objective function with the weight
function $w(\lambda )=\frac{1}{1+\lambda ^{2}}$ had been used in \cite%
{GAH2002} for the estimation of the Gaussian processes in stationary and
nonstationary cases respectively.
\end{rmrk}

\subsection{Minimum contrast estimation based on the Ibragimov contrast
function}

We consider now the minimum contrast functional motivated by the paper of
Ibragimov \cite{Ibragimov1967}, see also Anh, Leonenko and Sakhno \cite{ALS}.

We assume condition \textbf{A.I} \ and introduce the following condition

\noindent \textbf{B. I. }There exists a nonnegative function $w\left(
\lambda \right) ,\lambda $ $\in \mathbb{R},$ such that

(i) $w\left( \lambda \right) $ is symmetric about $\lambda =0:w\left(
\lambda \right) =w\left( -\lambda \right) ;$

(ii) $w\left( \lambda \right) f\left( \lambda ;\theta \right) $ is in $%
L_{1}\left( \mathbb{R}\right) $ for $\forall \theta \in \Theta .$

Under the condition B.I, we set
\begin{equation*}
\sigma ^{2}\left( \theta \right) =\int_{\mathbb{R}}f\left( \lambda ;\theta
\right) w\left( \lambda \right) d\lambda
\end{equation*}
and consider the factorization of the spectral density
\begin{equation*}
f\left( \lambda ;\theta \right) =\sigma ^{2}\left( \theta \right) \psi
\left( \lambda ;\theta \right) ,\quad \lambda \in \mathbb{R},\theta \in
\Theta .
\end{equation*}
For the function $\psi \left( \lambda ,\theta \right) ,\lambda $ $\in
\mathbb{R},$ $\theta \in \Theta ,$ we have
\begin{equation*}
\int_{\mathbb{R}}\psi \left( \lambda ;\theta \right) w\left( \lambda \right)
d\lambda =1
\end{equation*}
and we additionally suppose

\noindent \textbf{B. II.} The derivatives $\nabla _{{}\theta }\psi \left(
\lambda ;\theta \right) $ exist and
\begin{equation*}
\nabla _{{}\theta }\int_{\mathbb{R}}\psi \left( \lambda ;\theta \right)
w\left( \lambda \right) d\lambda =\int_{\mathbb{R}}\nabla _{{}\theta }\psi
\left( \lambda ;\theta \right) w\left( \lambda \right) d\lambda =0,
\end{equation*}
that is we can differentiate under the integral sign in the above integral.

Consider the following contrast process (or objective function):
\begin{equation}
U_{T}\left( \theta \right) =-\int_{\mathbb{R}}I_{T}\left( \lambda \right)
w\left( \lambda \right) \log \psi \left( \lambda ;\theta \right) d\lambda
,\quad \theta \in \Theta .  \label{(1i)}
\end{equation}
Define also the function
\begin{equation}
K\left( \theta _{0};\theta \right) =\int_{\mathbb{R}}f\left( \lambda ;\theta
_{0}\right) w\left( \lambda \right) \log \frac{\psi \left( \lambda ;\theta
_{0}\right) }{\psi \left( \lambda ;\theta \right) }d\lambda ,\quad \theta
_{0},\theta \in \Theta .  \label{(2i)}
\end{equation}

\noindent \textbf{B. III. }$f\left( \lambda ;\theta _{0}\right) w\left(
\lambda \right) \log \psi \left( \lambda ;\theta \right) \in L_{1}\left(
\mathbb{R}^{1}\right) \cap L_{2}\left( \mathbb{R}^{1}\right) ,\quad \forall
\theta \in \Theta .$

\noindent \textbf{B. IV.} There exists a function $v\left( \lambda \right)
,\lambda \in \mathbb{R}^{1},$ such that

(i) the function $h\left( \lambda ;\theta \right) =v\left( \lambda \right)
\log \psi \left( \lambda ;\theta \right) $ is uniformly continuous in $%
\mathbb{R}^{1}\times \Theta $;

(ii) $f\left( \lambda ;\theta _{0}\right) \frac{w\left( \lambda \right) }{%
v\left( \lambda \right) }\in L_{1}\left( \mathbb{R}^{1}\right) \cap
L_{2}\left( \mathbb{R}^{1}\right) .$

\begin{thrm}
\label{Ib1} Let conditions AI, B.I - B.IV be satisfied. Then the function $%
K\left( \theta _{0};\theta \right) $ defined by (\ref{(2i)}) is the contrast
function for the contrast process $U_{T}\left( \theta \right) $ defined by (%
\ref{(1i)}). Moreover the minimum contrast estimator $\widehat{\theta }_{T}$
defined as
\begin{equation}
\widehat{\theta }_{T}=\underset{\theta \in \Theta }{\arg \min }U_{T}\left(
\theta \right) ,  \label{(3i)}
\end{equation}%
is a consistent estimator of the parameter $\theta ,$ that is, $\widehat{%
\theta }_{T}\rightarrow \theta _{0}$ in $P_{0}$-probability as $T\rightarrow
\infty ,$ and the estimator
\begin{equation*}
\widehat{\sigma }_{T}^{2}=\int_{\mathbb{R}^{n}}I_{T}\left( \lambda \right)
w\left( \lambda \right) d\lambda
\end{equation*}%
is a consistent estimator of the parameter $\sigma ^{2}\left( \theta \right)
,$ that is, $\widehat{\sigma }_{T}^{2}\rightarrow \sigma ^{2}\left( \theta
_{0}\right) $ in $P_{0}$-probability as $T\rightarrow \infty .$
\end{thrm}

To formulate the result on the asymptotic distribution of the minimum
contrast estimator (\ref{(3i)}) we need some further conditions.

\noindent \textbf{B. V. }The function $\psi \left( \lambda ;\theta \right) $
is twice differentiable in a neighborhood of the point $\theta _{0}$ and

\noindent (i) $f\left( \lambda ;\theta \right) w\left( \lambda \right) \frac{%
\partial ^{2}}{\partial \theta _{i}\partial \theta _{j}}\log \psi \left(
\lambda ,\theta \right) $ $\in $ $L_{1}\left( \mathbb{R}\right) $ $\cap $ $%
L_{2}\left( \mathbb{R}\right) ,$ $i,j=1,...,m,$ $\theta \in \Theta ;$

\noindent (ii) $f\left( \lambda ;\theta _{0}\right) \in \mathbf{L}_{p}\left(
\mathbb{R}^{1}\right) ,$ $w\left( \lambda \right) \frac{\partial }{\partial
\theta _{i}}\log \psi \left( \lambda ,\theta \right) $ $\in $ $\mathbf{L}%
_{q}\left( \mathbb{R}^{1}\right) ,$\

for some $p,$ $q$ such that $\frac{1}{p}+\frac{1}{q}\leq \frac{1}{2}$, $%
i=1,...,m,$ $\theta \in \Theta ;$

\noindent (iii) $T^{1/2}\int_{\mathbb{R}^{1}}(EI_{T}\left( \lambda \right)
-f\left( \lambda ;\theta _{0}\right) )w\left( \lambda \right) \frac{\partial
}{\partial \theta _{i}}\log \psi \left( \lambda ;\theta \right) \ d\lambda
\rightarrow 0$ \ as $T\rightarrow \infty ,$

for all $i=1,...,m,\theta \in \Theta ;$

\noindent (iv) the second order derivatives $\frac{\partial ^{2}}{\partial
\theta _{i}\partial \theta _{j}}\log \psi \left( \lambda ;\theta \right)
,i=1,...,m,$ are continuous in $\theta .$

\noindent \textbf{B. VI. }The matrices $S\left( \theta \right) =\left(
s_{ij}\left( \theta \right) \right) _{i,j=1,...,m}$ and $A\left( \theta
\right) =\left( a_{ij}\left( \theta \right) \right) _{i,j=1,...,m}$ are
positive definite where
\begin{equation*}
s_{ij}\left( \theta \right) =\int_{\mathbb{R}}f\left( \lambda ;\theta
\right) w\left( \lambda \right) \frac{\partial ^{2}}{\partial \theta
_{i}\partial \theta _{j}}\log \psi \left( \lambda ;\theta \right) d\lambda
\end{equation*}
\begin{equation*}
=\sigma ^{2}\left( \theta \right) \int_{\mathbb{R}}w\left( \lambda \right) %
\left[ \frac{\partial ^{2}}{\partial \theta _{i}\partial \theta _{j}}\psi
\left( \lambda ,\theta \right) -\frac{1}{\psi \left( \lambda ,\theta \right)
}\frac{\partial }{\partial \theta _{i}}\psi \left( \lambda ,\theta \right)
\frac{\partial }{\partial \theta _{j}}\psi \left( \lambda ,\theta \right) %
\right] d\lambda ,
\end{equation*}
\begin{equation*}
a_{ij}\left( \theta \right) =4\pi \int_{\mathbb{R}^{1}}f^{2}\left( \lambda
;\theta \right) w^{2}\left( \lambda \right) \frac{\partial }{\partial \theta
_{i}}\log \psi \left( \lambda ;\theta \right) \frac{\partial }{\partial
\theta _{j}}\log \psi \left( \lambda ;\theta \right) d\lambda
\end{equation*}
\begin{equation*}
+2\pi \frac{d_{4}}{d_{2}^{2}}\int_{\mathbb{R}}\frac{w\left( \lambda \right)
f\left( \lambda ;\theta \right) }{\psi \left( \lambda ;\theta \right) }\frac{%
\partial }{\partial \theta _{i}}\psi \left( \lambda ;\theta \right) d\lambda
\int_{\mathbb{R}}\frac{w\left( \lambda \right) f\left( \lambda ;\theta
\right) }{\psi \left( \lambda ;\theta \right) }\frac{\partial }{\partial
\theta _{j}}\psi \left( \lambda ;\theta \right) d\lambda
\end{equation*}
\begin{eqnarray*}
&=&4\pi \left( \sigma ^{2}\left( \theta \right) \right) ^{2}\int_{\mathbb{R}%
}w^{2}\left( \lambda \right) \frac{\partial }{\partial \theta _{i}}\psi
\left( \lambda ;\theta \right) \frac{\partial }{\partial \theta _{j}}\psi
\left( \lambda ;\theta \right) d\lambda \\
&&+2\pi \frac{d_{4}}{d_{2}^{2}}\left( \sigma ^{2}\left( \theta \right)
\right) ^{2}\int_{\mathbb{R}}w\left( \lambda \right) \frac{\partial }{%
\partial \theta _{i}}\psi \left( \lambda ;\theta \right) d\lambda \int_{%
\mathbb{R}}w\left( \lambda \right) \frac{\partial }{\partial \theta _{j}}%
\psi \left( \lambda ;\theta \right) d\lambda
\end{eqnarray*}

\begin{thrm}
\label{Ib2} Let the conditions AI, B.I - B.VI be satisfied. Then as $%
T\longrightarrow \infty $%
\begin{equation*}
T^{1/2}\left( \widehat{\theta _{T}}-\theta _{0}\right) \overset{\mathcal{D}}{%
\rightarrow }N_{m}\left( 0,S^{-1}\left( \theta _{0}\right) A\left( \theta
_{0}\right) S^{-1}\left( \theta _{0}\right) \right) ,
\end{equation*}
where $N_{m}\left( \cdot ,\cdot \right) $ denotes the $m$-dimensional
Gaussian law.
\end{thrm}

\textbf{Proofs of Theorems \ref{Wh1}, \ref{Wh2}, \ref{Ib1}, \ref{Ib2}.}

The results on consistency of estimators \ (Theorems \ref{Wh1} and \ref{Ib1}
) are consequences of corresponding theorems stated for the general case in
\cite{LS2006} for the Whittle functional and in \cite{ALS}, \cite{ALS1} for
the case of Ibragimov functional. We present here reasonings for the proofs
of Theorems \ref{Wh2} and \ref{Ib2}, which make use of CLT for bilinear
forms (Theorem \ref{linbilin} above). For the proofs the standard arguments
based on Taylor's formula for $\nabla _{{}\theta }U_{T}\left( \widehat{%
\theta }_{T}\right) $ are used. Namely, we can write the relation
\begin{equation*}
\nabla _{{}\theta }U_{T}\left( \widehat{\theta }_{T}\right) =\nabla
_{{}\theta }U_{T}\left( \theta _{0}\right) +\nabla _{{}\theta }\nabla
_{{}\theta }^{\prime }U_{T}\left( \theta _{T}^{\ast }\right) \left( \widehat{%
\theta }_{T}-\theta _{0}\right) ,
\end{equation*}
where $\left| \theta _{T}^{\ast }-\theta _{0}\right| <\left| \widehat{\theta
}_{T}-\theta _{0}\right| .$

It follows from the definition of minimum contrast estimators that for
sufficiently large $T$
\begin{equation*}
\nabla _{{}\theta }U_{T}\left( \theta _{0}\right) =-\nabla _{{}\theta
}\nabla _{{}\theta }^{\prime }U_{T}\left( \theta _{T}^{\ast }\right) \left(
\widehat{\theta }_{T}-\theta _{0}\right) ,
\end{equation*}
therefore, to state the asymptotic normality for the estimator $\widehat{%
\theta }_{T}$ , by Slutsky's arguments, one needs to deduce: (1) limit in
probability for $\nabla _{{}\theta }\nabla _{{}\theta }^{\prime }U_{T}\left(
\theta _{T}^{\ast }\right) $ and (2) limiting normal law for $T^{1/2}\nabla
_{{}\theta }U_{T}\left( \theta _{0}\right) .$

For the 1-st task we can use the same arguments as in the mentioned above
papers, and to rewrite (simplify) corresponding conditions for the case of
linear processes.

However, for the step (2) we can appeal now to Theorem \ref{linbilin}. We
provide the details below.

Consider firstly the case of Whittle functional. Limit in $P_{0}$%
-probability for $\nabla _{{}\theta }\nabla _{{}\theta }^{\prime
}U_{T}\left( \theta _{T}^{\ast }\right) $ is given by the matrix $%
W_{1}\left( \theta _{0}\right) .$

Next, consider
\begin{equation*}
\nabla _{{}\theta }U_{T}\left( \theta _{0}\right) =\frac{1}{4\pi }\int_{%
\mathbb{R}}\left( \left. \nabla _{{}\theta }\log f\left( \lambda ;\theta
\right) \right| _{\theta =\theta _{0}}+\left. \nabla _{{}\theta }\left(
\frac{1}{f\left( \lambda ;\theta \right) }\right) \right| _{\theta =\theta
_{0}}I_{T}\left( \lambda \right) \right) w\left( \lambda \right) d\lambda ,
\end{equation*}
which can be written in the form
\begin{eqnarray*}
\nabla _{{}\theta }U_{T}\left( \theta _{0}\right) &=&\left( J_{T}(\varphi
_{i})-J(\varphi _{i})\right) _{i=1,...,m} \\
&=&\left( \int_{\mathbb{R}}\varphi _{i}\left( \lambda ;\theta _{0}\right)
I_{T}\left( \lambda \right) d\lambda -\int_{\mathbb{R}}\varphi _{i}\left(
\lambda ;\theta _{0}\right) f\left( \lambda ;\theta _{0}\right) d\lambda
\right) _{i=1,...,m}
\end{eqnarray*}
where
\begin{equation*}
\varphi _{i}=\varphi _{i}\left( \lambda ;\theta _{0}\right) =-\text{ }\frac{1%
}{4\pi }\frac{1}{f^{2}\left( \lambda ;\theta _{0}\right) }w\left( \lambda
\right) \left. \left( \frac{\partial }{\partial \theta _{i}}f\left( \lambda
;\theta \right) \right) \right| _{\theta =\theta _{0}},\quad i=1,...,m.
\end{equation*}
Under the assumptions of Theorem \ref{Wh2} (see A.IV(ii)) in view of Theorem %
\ref{linbilin} we have the convergence
\begin{equation}
T^{1/2}\left( J_{T}(\varphi _{i})-EJ_{T}(\varphi _{i})\right) _{i=1,...,m}%
\overset{\mathcal{D}}{\rightarrow }N_{m}\left( 0,W_{2}\left( \theta
_{0}\right) +V\left( \theta _{0}\right) \right) ,  \label{11}
\end{equation}
where the the matrices $W_{2}\left( \theta _{0}\right) $ and $V\left( \theta
_{0}\right) $ are defined in the assumption A.V.

Further, in view of the assumption A.IV(iii)
\begin{equation*}
T^{1/2}\left( EJ_{T}(\varphi _{i})-J(\varphi _{i})\right) \overset{}{%
\rightarrow }0,\text{ as }T\rightarrow \infty
\end{equation*}
which, combined with (\ref{11}), implies
\begin{equation*}
T^{1/2}\left( J_{T}(\varphi _{i})-J(\varphi _{i})\right)
_{i=1,...,m}=T^{1/2}\nabla _{{}\theta }U_{T}\left( \theta _{0}\right)
\overset{\mathcal{D}}{\rightarrow }N_{m}\left( 0,W_{2}\left( \theta
_{0}\right) +V\left( \theta _{0}\right) \right) .
\end{equation*}

The case of Ibragimov functional is treated analogously. We have that \ $%
\nabla _{{}\theta }\nabla _{{}\theta }^{\prime }U_{T}\left( \theta
_{T}^{\ast }\right) $ converges \ in $P_{0}$-probability to the matrix $%
S\left( \theta _{0}\right) .$ Further,
\begin{equation*}
\nabla _{{}\theta }U_{T}\left( \theta _{0}\right) =-\int_{\mathbb{R}%
}I_{T}\left( \lambda \right) \left. \nabla _{{}\theta }\log \psi \left(
\lambda ;\theta \right) \right| _{\theta =\theta _{0}}w\left( \lambda
\right) d\lambda .
\end{equation*}
In view of B.II
\begin{equation*}
\int_{\mathbb{R}}f\left( \lambda ;\theta _{0}\right) \left. \nabla
_{{}\theta }\log \psi \left( \lambda ;\theta \right) \right| _{\theta
=\theta _{0}}w\left( \lambda \right) d\lambda =0,
\end{equation*}
and we can write
\begin{eqnarray*}
\nabla _{{}\theta }U_{T}\left( \theta _{0}\right) &=&\left( J_{T}(\varphi
_{i})-J(\varphi _{i})\right) _{i=1,...,m} \\
&=&\left( \int_{\mathbb{R}}\varphi _{i}\left( \lambda ;\theta _{0}\right)
I_{T}\left( \lambda \right) d\lambda -\int_{\mathbb{R}}\varphi _{i}\left(
\lambda ;\theta _{0}\right) f\left( \lambda ;\theta _{0}\right) d\lambda
\right) _{i=1,...,m},
\end{eqnarray*}
where now
\begin{equation*}
\varphi _{i}=\varphi _{i}\left( \lambda ;\theta _{0}\right) =w\left( \lambda
\right) \left. \frac{\partial }{\partial \theta _{i}}\log \psi \left(
\lambda ;\theta \right) \right| _{\theta =\theta _{0}},\quad i=1,...,m.
\end{equation*}
Again in view of Theorem \ref{linbilin}, under the assumption B.V(ii), we
obtain the convergence
\begin{equation*}
T^{1/2}\left( J(\varphi _{i})-EJ_{T}(\varphi _{i})\right) \overset{\mathcal{D%
}}{\rightarrow }N_{m}\left( 0,A\left( \theta _{0}\right) \right) ,\text{ as }%
T\rightarrow \infty ,
\end{equation*}
where the matrix $A\left( \theta _{0}\right) $ is defined in B.VI. By
assumption B.V(iii) the convergence $T^{1/2}\nabla _{{}\theta }U_{T}\left(
\theta _{0}\right) \overset{\mathcal{D}}{\rightarrow }N_{m}\left( 0,A\left(
\theta _{0}\right) \right) $ follows.

\appendix

\section{Fractional Riesz-Bessel Motion}

In this Appendix we mainly review a number results discussed in Gay and
Heyde \cite{GayHeyde1990}, Anh, Angulo and Ruiz-Medina \cite{AARM1999}, Anh,
Leonenko and Mc Vinish \cite{ALMcV2001}, Anh and Leonenko \cite%
{AnhLeonenko2001}, \cite{AnhLeonenko2002}, Kelbert, Leonenko and Ruiz-Medina
\cite{KLR} (see also references therein). Also we introduce a not
necessarily Gaussian Riesz-Bessel stationary process and formulate the
central limit theorem for such a processes as well as for quadratic forms of
such a processes.

The fractional operators are natural mathematical objects to describe the
singular phenomena of random fields such as long range dependence or/and
intermittency.

In particular Gay and Heyde \cite{GayHeyde1990} introduced a class of random
fields as solutions of fractional Helmholtz equation driven by white noise,
contained the fractional operator $(cI-\Delta )^{\alpha /2},\;c\geq 0$ (and
its limit as $c\rightarrow 0\quad (-\Delta )^{\alpha /2})$, where $\Delta $
is the $d$-dimensional Laplacian and $I$ is the identity operator (see also
\cite{KLR} for properties of such fields and possible generalization). Anh,
Angulo, Ruiz-Medina \cite{AARM1999} (see also \cite{ALMcV2001}, \cite%
{AnhLeonenko2001}, \cite{AnhLeonenko2002} and references therein)
generalized the fractional stochastic equation of Gay and Heyde in which the
fractional Helmgoltz operator $(cI-\Delta )^{\alpha /2},\;c\geq 0$ or the $d$%
-dimensional Laplacian ($c\rightarrow 0$) is replaced by a fractional
Laplace-type operator of the form $-\left( I-\Delta \right) ^{\gamma
/2}(-\Delta )^{\alpha /2},\alpha >0,\;\gamma \geq 0,$ where the operators $%
-\left( I-\Delta \right) ^{\gamma /2},\gamma \geq 0,$ and $(-\Delta
)^{\alpha /2},\;\alpha >0,$ are interpreted as inverses to the Bessel and
Riesz potentials (see \cite{Stein1970}, pp.~134-138), that is integral
operators, whose kernels have a Fourier transforms $(2\pi )^{-d/2}\left(
1+\left\Vert \lambda \right\Vert ^{2}\right) ^{-\gamma /2},\;\lambda \in
\mathbb{R}^{d},$ and $(2\pi )^{-d/2}\left\Vert \lambda \right\Vert ^{-\alpha
},\;\lambda \in \mathbb{R}^{d},$ respectively. Then there exists a
generalized random field $\zeta (x),\;x\in \mathbb{R}^{d},$ on fractional
Sobolev space, which is defined by the equation
\begin{equation}
\left( I-\Delta \right) ^{\gamma /2}(-\Delta )^{\alpha /2}\zeta
(x)=e(x),\;x\in \mathbb{R}^{d},  \label{ap1}
\end{equation}%
where $\left\{ e(x),\;x\in \mathbb{R}^{d}\right\} $ is a Gaussian white
noise or equivalently (in the sense of second-order moments) there exists a
random field with the spectral density
\begin{equation}
f(\lambda )=\frac{c}{\left\Vert \lambda \right\Vert ^{2\alpha }\left(
1+\left\Vert \lambda \right\Vert ^{2}\right) ^{\gamma }},\;\lambda \in
\mathbb{R}^{d},\;c>0.  \label{ap2}
\end{equation}

These random fields were named the fractional Riesz-Bessel motion.

For the random fields with stationary increments we assume $\alpha \in
\left( \frac{d}{2},\frac{d}{2}+1\right) $, $\gamma \geq 0.$ In particular,
for $d=1,$ there exists a Gaussian stochastic process with stationary
increments and the spectral density (\ref{ap2}), where $\alpha \in \left(
\frac{1}{2},\frac{3}{2}\right) ,\;\gamma \geq 0.$ This fractional
Riesz-Bessel motion (FRBM) is a generalization of the fractional Brownian
motion (FBM) (see, for instance, Samorodnitsky and Taqqu \cite{ST1994}). FBM
is a limiting case of the Riesz-Bessel (non-stationary) motion with $\gamma
=0$ (in terms of the Hurst parameter $H\in (0,1),$ the spectral density of
the FBM with long-range dependence ($H\in (\frac{1}{2},0)$) is equal to $%
\frac{1}{\left\vert \lambda \right\vert ^{2H+1}}).$ The FRBM is not
self-similar (unless when $\gamma =0)$, but it is locally self-similar.

For $d\geq 1$ the presence of the Bessel operator is essential for a study
of homogeneous (and isotropic) solutions of (\ref{ap1}) with spectral
density (\ref{ap2}), which requires $0\leq \alpha <\frac{d}{2},\;\alpha
+\gamma >\frac{d}{2};$ that is the condition $\gamma >0$ is necessary for $%
f(\lambda )\in L_{1}(\mathbb{R}^{d}).$ Thus the homogeneous isotropic FRBM
can be defined as a Gaussian random field with zero mean and covariance
function of the form
\begin{equation}
B_{\alpha ,\gamma }(x)=\int\limits_{\mathbb{R}^{d}}e^{i<\lambda ,x>}\frac{c}{%
\left\Vert \lambda \right\Vert ^{2\alpha }\left( 1+\left\Vert \lambda
\right\Vert ^{2}\right) ^{\gamma }}d\lambda ,\;x\in \mathbb{R}^{d},
\label{ap3}
\end{equation}%
where $0\leq \alpha <\frac{d}{2},\;\alpha +\gamma >\frac{d}{2}.$ Note that
for $\alpha =0$, the covariance structure (\ref{ap3}) belongs the Mat\'{e}rn
class, that is with
\begin{equation*}
c=\frac{\Gamma (\gamma )}{\pi ^{d/2}2^{d-1}\Gamma \left( \frac{2\gamma -d}{2}%
\right) };
\end{equation*}%
\begin{equation}
B_{\alpha ,\gamma }(x)=\frac{1}{2^{\frac{2-d}{2}}\Gamma \left( \frac{2\gamma
-d}{2}\right) }\frac{K_{\frac{2\gamma -d}{2}}}{\left\Vert x\right\Vert ^{%
\frac{d-2\gamma }{2}}},\;x\in \mathbb{R}^{d},\;\gamma >\frac{d}{2},
\label{ap4}
\end{equation}%
where
\begin{equation*}
K_{\nu }(z)=\frac{1}{2}\int\limits_{0}^{\infty }s^{\nu -1}\exp \left\{ -%
\frac{1}{2}\left( s+\frac{1}{s}\right) z\right\} ds,\;z\geq 0,\;\nu \in R,
\end{equation*}%
is the modified Bessel function of the third kind of order $\nu $ or Mc
Donald's function. Note that
\begin{equation*}
K_{\nu }(z)=K_{-\nu }(z),\quad K_{\nu }(z)\sim \Gamma (\nu )2^{\nu
-1}z^{-\nu },
\end{equation*}%
for$\;\nu >0$ as $z\rightarrow 0,$%
\begin{equation*}
K_{\frac{1}{2}}(z)=\sqrt{\frac{\pi }{2}}\frac{e^{-z}}{\sqrt{z}}.
\end{equation*}

Thus, we have $B_{\alpha ,\gamma }(0)=1.$

Note that for $d=1,\;\alpha =0,\;\gamma =1,$ the covariance structure (\ref%
{ap4}) becomes $B_{0,1}(x)=e^{-x},\;x\geq 0,$ that is stationary Gaussian
Riesz-Bessel motion is identical to the Gaussian Ornstein-Uhlenbeck process.

\begin{rmrk}
These results can be generalized to the case when the above fractional
operator is replaced by the operator
\begin{equation*}
H=\frac{\partial ^{\beta }}{\partial t^{\beta }}+\mu \left( I-\Delta \right)
^{\gamma /2}(-\Delta )^{\alpha /2},\;0\leq \beta \leq 2,\ \alpha >0,\;\gamma
\geq 0,
\end{equation*}%
where $\frac{\partial ^{\beta }}{\partial t^{\beta }}$ is the regularized
fractional derivative. In particular, the Green function of the fractional
heat equation:$\;Hu(t,x)=0,$ $t>0,\;x\in \mathbb{R}^{d},$ can be given as
inverse Fourier transform of the function
\begin{equation*}
E_{\beta }\left( -\mu t^{\beta }\left\Vert \lambda \right\Vert ^{\alpha
}(1+\left\Vert \lambda \right\Vert ^{2})^{\gamma /2}\right) ,\;t>0,\;\lambda
\in \mathbb{R}^{d},
\end{equation*}%
where
\begin{equation*}
E_{\beta }(z)=\sum_{j=1}^{\infty }\frac{z^{j}}{\Gamma (\beta _{j}+1)},\;z\in
\text{\texttt{C}},\;\beta >0
\end{equation*}%
is the Mittag-Leffler function (see \cite{AnhLeonenko2001}, \cite%
{AnhLeonenko2002} for details and references).
\end{rmrk}

In order to introduce a Riesz-Bessel motion driven by L\'{e}vy noise, we
restrict our attention to the stationary case and $d=1$ (replacing the space
parameter $x$ into $t)$. For the function
\begin{equation*}
a(\lambda )=\frac{\sqrt{c}}{(i\lambda )^{\alpha }(1+i\lambda )^{\gamma }}%
,\;\lambda \in \mathbb{R},\;a+\gamma >1,\;\alpha \geq 0,
\end{equation*}
we introduce the function
\begin{equation}
\hat{a}(t)=\int_{\mathbb{R}}e^{it\lambda }a(\lambda )d\lambda =\left\{
\begin{array}{ll}
\frac{2\pi }{\Gamma (\alpha +\gamma )}t^{\alpha +\gamma
-1}e^{-1}\,_{1}F_{1}(\gamma ,\alpha +\gamma ;t), & t\geq 0, \\
0, & t<0,%
\end{array}
\right. ,  \label{ap5}
\end{equation}
where the confluent hypergeometric function
\begin{equation*}
_{1}F_{1}(a,b;z)=\sum_{n=0}^{\infty }\frac{(a)_{n}}{(b)_{n}}\frac{z^{n}}{n!}%
,\;(c)_{n}=c(c+1)\cdots (c+n-1),\;(c)_{0}=1.
\end{equation*}

The Riesz-Bessel motion driven by L\'{e}vy noise can be defined as the
linear process
\begin{equation}
X_{t}=\int_{\mathbb{R}}\hat{a}(t-s)d\xi (s),  \label{ap6}
\end{equation}
where $\xi (t),\;t\in R$ is a L\'{e}vy process with cumulant function
\begin{equation*}
\kappa (z)=\log E\exp \left\{ iz\xi (1)\right\} ,
\end{equation*}
such that $\kappa ^{(k)}(0)\neq 0,\;k\geq 2,$ and $\hat{a}(.)$ is defined by
(\ref{ap5}). The $k$-th order spectral densities of the Riesz-Bessel motion
driven by L\'{e}vy noise (\ref{ap6}) take the form:
\begin{equation}
f_{k}(\lambda _{1},\ldots ,\lambda _{k-1}) =(2\pi )^{-k+1}i^{-k}\kappa
^{(k)}(0)a(\lambda _{1})\cdots a(\lambda _{k-1})\overline{a(\lambda
_{1}+\cdots +\lambda _{k-1})},  \label{ap7}
\end{equation}
which reduces to the second-order spectral density
\begin{eqnarray}
f_{2}(\lambda ) &=&\frac{c}{\left| \lambda \right| ^{2\alpha }(1+\lambda
^{2})^{\gamma }},  \label{ap8} \\
c &=&\frac{\kappa ^{(2)}(0)}{2\pi },\;0\leq \alpha <\frac{1}{2},\;\alpha
+\gamma >\frac{1}{2},\;\lambda \in \mathbb{R}^{1}.  \notag
\end{eqnarray}

For the Gaussian case, of course, $\kappa ^{(k)}(0)=0,\;k\geq 3.$

Note that for $\alpha =0,\;\gamma =1$ we arrive to the Ornstein-Uhlenbeck
process driven by L\'{e}vy noise (\cite{AHL}).

As a consequence of the Theorems of Sections 4 and 5 we obtain the following
result for the linear process (\ref{ap6}). (Cf. also with Theorem 4.1)

\begin{thrm}
Consider the Riesz-Bessel stationary motion (\ref{ap6}) and assume that all
cumulants of L\'{e}vy process are finite. Let
\begin{equation*}
S_{T}=\int_{-T/2}^{T/2}X_{s}ds,\text{ \ \ \ }Q_{T}=\int_{-T/2}^{T/2}%
\int_{-T/2}^{T/2}\hat{b}(t-s)\left[ X_{t}X_{s}-\mathit{E}X_{t}X_{s}\right]
dtds,
\end{equation*}

then:

i) if
\begin{equation*}
\alpha +\gamma >\frac{1}{2},\text{ \ }\alpha \leq 0,
\end{equation*}

then the central limit theorem holds:
\begin{equation*}
T^{-1/2}S_{T}\rightarrow N(0,\sigma ^{2}),T\rightarrow \infty ,
\end{equation*}

where
\begin{equation*}
\sigma ^{2}=\kappa ^{(2)}(0);
\end{equation*}

\bigskip

ii) if for some $p>1,q>1,$ such that $\frac{1}{p}+\frac{1}{q}\leq \frac{1}{2}%
,$ we have
\begin{equation*}
b(\lambda )\in L_{q},\alpha +\gamma >\frac{1}{2p},\text{ \ }\alpha <\frac{1}{%
2p},
\end{equation*}

then the central limit theorem holds:
\begin{equation*}
T^{-1/2}Q_{T}\rightarrow N(0,\sigma ^{2}),T\rightarrow \infty ,
\end{equation*}

where
\begin{equation*}
\sigma ^{2} =2\kappa ^{(2)}(0)\int_{\mathbb{R}}b^{2}(\lambda )\frac{1}{%
\left\vert \lambda \right\vert ^{4\alpha }(1+\lambda ^{2})^{2\gamma }}%
d\lambda +\kappa ^{(4)}(0)\frac{1}{(2\pi )^{2}}\int_{\mathbb{R}}b(\lambda )%
\frac{1}{\left\vert \lambda \right\vert ^{2\alpha }\left( 1+\lambda
^{2}\right) ^{\gamma }}d\lambda .
\end{equation*}
\end{thrm}

\section{Kernel estimates \label{s:ker}}

Consider the Dirichlet type kernel
\begin{equation*}
\Delta _{T}(\lambda )=\int_{t\in I_{T}}e^{it\lambda }\nu (dt).
\end{equation*}

When $d=1$ and $I_{1}=[-1/2,1/2]$ ($I_T=I_1 T$), one gets the clasical
discrete/continuous time Dirichlet kernels:
\begin{equation*}
\Delta _{T}(\lambda )=\sum_{-T/2}^{T/2}e^{it\lambda }=\frac{\sin
((T+1)\lambda /2)}{\sin (\lambda /2)},\quad \Delta _{T}(\lambda
)=\int_{-T/2}^{T/2}e^{it\lambda }dt=\frac{\sin (T\lambda /2)}{\lambda /2},
\end{equation*}%
respectively. For general $d$ and $I_{T}=[-T/2,T/2]^{d}$, putting $\lambda
=(\lambda _{1},..,\lambda _{d})$, it follows that $\Delta _{T}(\lambda
)=\prod_{j=1}^{d}\Delta _{T}(\lambda _{j})$.

Note that in the continuous case, by scaling, one finds
\begin{equation}
|\!|\Delta _{T}(\lambda )|\!|_{p}=T^{1-1/p}C_{p},\quad \;\lambda \in \mathbb{%
R},\quad 1<p<\infty ,
\end{equation}%
\begin{equation}
|\!|\Delta _{T}(\lambda )|\!|_{p}=T^{d(1-1/p)}C_{p}^{d},\quad \lambda \in
\mathbb{R}^{d},\quad 1<p<\infty.
\end{equation}
with $C_{p}=(2\int_{R}|\frac{\sin (z)}{z}|^{p}dz)^{\frac{1}{p}}.$

In the discrete case, similar estimates may be obtained by using the
inequality%
\begin{equation*}
\left\vert \frac{\sin ((T+1)\lambda /2)}{\sin (\lambda /2)}\right\vert \leq
\tilde{C}\frac{T}{1+T\left\vert \lambda \right\vert },\lambda \in \lbrack
-\pi ,\pi ).
\end{equation*}

We find then:

\begin{equation*}
|\!|\Delta _{T}(\lambda )|\!|_{p}\leq T_{p}^{1-1/p}\tilde{C}^{\frac{1}{p}%
}\left( \int\limits_{\mathbb{R}}\frac{dz}{\left( 1+\left\vert z\right\vert
^{p}\right) }\right) ^{\frac{1}{p}}\;,\ \ \lambda \in \lbrack -\pi ,\pi
),1<p<\infty ,\quad
\end{equation*}%
\begin{equation*}
|\!|\Delta _{T}(\lambda )|\!|_{p}\leq T_{p}^{d(1-1/p)}\left[ \tilde{C}^{%
\frac{1}{p}}\left( \int\limits_{\mathbb{R}}\frac{dz}{\left( 1+\left\vert
z\right\vert ^{p}\right) }\right) ^{\frac{1}{p}}\right] ^{d}\;,\ \ \lambda
=(\lambda _{1},..,\lambda _{d})\in \lbrack -\pi ,\pi )^{d},1<p<\infty.
\end{equation*}

In the case of the Euclidean ball $I_{T}=B_{T}=\{t\in \mathbb{R}%
^{d}:\left\Vert t\right\Vert \leq T/2\},$ we find again by scaling in the
continuous case%
\begin{equation*}
\Delta _{T}(\lambda )=\int_{B_{T}}e^{it\lambda }dt=(2\pi )^{\frac{d}{2}%
}J_{d/2}(\left\Vert \lambda \right\Vert \frac{T}{2})/\left\Vert \lambda
\right\Vert ^{d/2},\quad \lambda \in \mathbb{R}^{d},
\end{equation*}%
where $J_{\nu }(z)$ is the Bessel function of the first kind and order $\nu
. $ It is known that $J_{\nu }(z)\leq const/\sqrt{z}$ for a large $z$, thus
for the ball%
\begin{equation*}
\left\Vert \Delta _{T}(\lambda )\right\Vert _{p}=\left\{
\begin{array}{ccc}
T^{(1-\frac{1}{p})}C_{p}, & d=1, & p>1, \\
T^{d(\frac{1}{2}-\frac{1}{p})}\bar{C}_{p}, & d\geq 2, & p>\frac{2d}{d+1},%
\end{array}%
\right.
\end{equation*}%
\begin{equation*}
\bar{C}_{p}=2^{^{d(\frac{1}{2}-\frac{1}{p})}}\left( 2\pi \right) ^{\frac{d}{2%
}}\left\vert s(1)\right\vert \left( \int_{0}^{\infty }\rho ^{d-1}\left\vert
\frac{J_{\frac{d}{2}}(\rho )}{\rho ^{d/2}}\right\vert ^{p}\right) ^{1/p},
\end{equation*}%
where $\left\vert s(1)\right\vert $ is the surface area of the unit ball in $%
\mathbb{R}^{d},~d\geq 2.$

Similar estimates may be obtained for the $L_{p\text{ }}(1<p\leq \infty )$
norms of the discrete Dirichlet kernel.

\textbf{Note:} These results are particular cases of the so-called
Hardy-Littlewood Theorem (see, for instance, Zigmund \cite{Z}, V. II, XII,
\S 6), which can be formulated as follows:

\begin{thrm}
\label{t1.4 copy(2)} Let $a_{n}\geq a_{n+1}\geq \cdots $ and $%
a_{n}\rightarrow 0$. Consider the series
\begin{equation}
\sum\limits_{n=1}^{\infty }a_{n}\cos n\lambda  \label{(1)}
\end{equation}%
and
\begin{equation}
\sum\limits_{n=1}^{\infty }a_{n}\sin n\lambda  \label{(2)}
\end{equation}%
and define by $f(\lambda)$ and $g(\lambda)$ the sums of the series (\ref{(1)}%
) and (\ref{(2)}) respectively at the points where the series converge. A
necessary and sufficient condition that the function $f$ (or $g$) belongs to
$L_{p}$, $1<p<\infty ,$ is the following
\begin{equation*}
\sum\limits_{n=1}^{\infty }a_{n}^{p}n^{p-2}\,\,<\,\infty .
\end{equation*}%
Moreover,
\begin{equation*}
\Vert f\Vert _{p}^{p}\asymp \sum\limits_{n=1}^{\infty }a_{n}^{p}n^{p-2}.
\end{equation*}
\end{thrm}

Clearly, $1=1=\cdots =1>0=\cdots $ is nonincreasing, and we arrive thus to
the following estimate for Dirichlet kernels:

\begin{equation*}
|\!|\sum_{t=1}^{T}e^{it\lambda }|\!|_{p}\!\leq C\left(
\sum_{t=1}^{T}1^{p}t^{p-2}\right) ^{\frac{1}{p}}\leq CT^{\frac{p-1}{p}%
},1<p<\infty ,
\end{equation*}

and
\begin{equation*}
\int_{0}^{1}\left\vert \sum_{t=0}^{T-1}e^{2\pi it\lambda }\right\vert
^{p}d\lambda =T^{p-1}\frac{2}{\pi }\int_{0}^{\infty }\left\vert \frac{\sin u%
}{u}\right\vert ^{p}du+R_{p}(T),
\end{equation*}%
where the error term%
\begin{equation*}
R_{p}(T)=\left\{
\begin{array}{ll}
O_{p}(T^{p-3}), & p>3 \\
O(\log T), & p=3 \\
O_{p}(1), & 1<p<3,%
\end{array}%
\right.
\end{equation*}%
where $O_{p}$ means that constants depend on $p.$(See, e.g., \cite{anderson}%
.)

Note that for $p=1$

\begin{equation*}
\int_{0}^{1}\left\vert \sum_{t=0}^{T-1}e^{2\pi it\lambda }\right\vert
d\lambda \asymp \frac{4}{\pi ^{2}}\log T.
\end{equation*}

\section{The homogeneous H\"older-Young-Brascamp-Lieb inequality \label%
{s:HYBu}}

Subtle modifications of the conditions of the H\"{o}lder inequality must be
made when the arguments of the functions involved are restricted to some
subspaces \cite{Fri}. Starting with Brascamp and Lieb \cite{BL} and Lieb
\cite{L} (who considered only the case $S=\mathbb{R}$), and following with
Ball \cite{Bal}, Barthe \cite{Barthe} and Carlen, Loss and Lieb \cite{CLL},
this generalization of the classical inequalities of H\"{o}lder and Young
seems to have attained now its definite form in the work of Bennett,
Carbery, Christ and Tao \cite{BCCT1}, \cite{BCCT}.

We review now a particular case of this result.

Let
\begin{equation*}
{x}=(x_{1},\ldots ,x_{m})\in S^{m}
\end{equation*}%
where $S$ may be either the multidimensional torus, integers or reals
\begin{equation*}
S=%
\begin{cases}
\lbrack -\pi ,\pi ]^{d} \\
\mathbb{Z}^{d} \\
\mathbb{R}^{d}%
\end{cases}%
\end{equation*}%
endowed with the respective normalized Haar measure $\mu (dx)$.

When $d=1,$ the convergence of integrals of the form:
\begin{equation*}
\boxed{\int_{x \in S^m} \frac{\mu(dx_1) ...\mu(dx_m)}
{l_1(x)^{z_1}...l_k(x)^{z_k}}}
\end{equation*}%
where $(l_{1},...,l_{k})$ are linear transformations
\begin{equation*}
l_{j}:S^{m}\;->\;S,\quad l_{j}(x)=<\alpha _{j},x>,\quad j=1,...,k
\end{equation*}%
and where in the first two cases $\alpha _{j}$ are supposed to have integer
coefficients, is a fundamental question arising in many applications.

Let $M$ denote the matrix with columns $\alpha _j, j=1,...,k$. It was for
long known to physicists that, when $M$ is fixed, convergence holds for $%
z=(z_1,...,z_k)$ belonging to a certain ``power counting polytope" PCP
(these are relatively similar in all the three cases -- see Theorem \ref%
{t:HYB} below).

It was first noticed in \cite{AB} and \cite{AT}, in the easier case of
\textbf{unimodal matrices} $M$, that under the same \textquotedblleft power
counting conditions" on $z_{j}=p_{j}^{-1}$, a H\"{o}lder-type inequality
\begin{equation*}
\biggl|\int_{S^{m}}\prod_{j=1}^{k}f_{j}(l_{j}({x}))\mu (d{x})\biggl|\leq
K\prod_{j=1}^{k}\Vert f_{j}\Vert _{p_{j}}\leqno(GH)
\end{equation*}%
holds, with the powers in (\ref{conv}) being replaced by arbitrary functions
satisfying integrability conditions $f_{j}\in L_{p_{j}},j=1,...,k$, and with
$K=1$. Note that Brascamp and Lieb \cite{BL} had already studied the analog
harder inequality for general matrices (in the case $S=\mathbb{R}^{d}$), but
without pinpointing exactly the polytope; this was done later by Barthe \cite%
{Barthe}.

For an example, consider the integral
\begin{equation*}
J= \int_S \int_S f_1(x_1)f_2(x_2)f_3(x_1+x_2)f_4(x_1-x_2) dx_1dx_2
\end{equation*}
where $S=\mathbb{R}$. Here $m=2,\ k=4$ and the matrix
\begin{equation*}
M=
\begin{pmatrix}
1 & 0 & 1 & 1 \\
0 & 1 & 1 & -1%
\end{pmatrix}%
\end{equation*}
has rank $r(M)=2$. The theorem below will ensure that
\begin{equation*}
|J|\leq\Vert f_1\Vert_{1/z_1}\ \Vert f_2\Vert_{1/z_2}\ \Vert
f_3\Vert_{1/z_3} \Vert f_4\Vert_{1/z_4}
\end{equation*}
for any ${\ z}=(z_1,z_2,z_3,z_4) \in [0,1]^4$ satisfying $z_1+z_2+z_3+z_4 =
2 $, e.g. if ${\ z}=(0,1,1/2, 1/2)$, then
\begin{equation*}
|J|\leq\left( \sup_{0\leq x\leq 1}|f_1(x)|\right)\left( \int_0^1|f_2(x)|dx
\right) \left( \int_0^1 f_3^2(x)dx\right)^{1/2}\left( \int_0^1 f_4^2(x)dx
\right)^{1/2}.
\end{equation*}

It is easy to check (and true in general) that the extremal points of the
PCP have only $0$ and $1$ coordinates (which may be exploited for
establishing the result). Note also that the matrix $M$ in this example is
not unimodal; as a consequence, the optimal constant $K=K(z)$ is not $1$ at
all the extremal points, the exception being $(0,0,1,1)$, where it is $%
2^{-1} $; also, the functions achieving equality must be Gaussian (which
holds in general, cf. Brascamp-Lieb \cite{BL}).

We will formulate now simultaneously the H\"older-Young-Brascamp-Lieb
inequality in the three cases:

\begin{itemize}
\item[(C1)] $\quad \mu(d x_j)$ is normalized Lebesgue measure on the torus $%
[-\pi,\pi]^{d}$, and $M$ has all its coefficients integers.

\item[(C2)] $\quad \mu(d x_j)$ is counting measure on $\mathbb{Z}^{d}$, $M $
has all its coefficients integers, and is unimodular, i.e. all its
non-singular minors of dimension $m \times m$ have determinant $\pm 1$.

\item[(C3)] $\quad \mu(d x_j)$ is Lebesgue measure on $(-\infty, +\infty
)^{d}$.
\end{itemize}

The result below specifies the domain of validity of H\"older's inequality,
called \textbf{power counting polytope}, in terms of linear inequalities
involving the rank $r(A)$ of arbitrary subsets $A$ of columns of the matrix $%
M$ (including the empty set $\emptyset $). It is also possible to express
the inequalities in terms of the ``dual rank" $r^*(A)$ defined by a dual
matrix $M^*$ whose lines are orthogonal to those of $M$, by using the
\textbf{duality relation}
\begin{equation*}
r^*(A)=|A|-r(M) +r(A^c),\ \forall A
\end{equation*}

\begin{thrm}
\label{t:HYB} \textbf{(Homogeneous H\"older-Young-Brascamp-Lieb inequality)}%
. Let $l_j(x)=x^t \alpha _j, j=1, ...,k$ be linear functionals $l_j:S^m\; ->
S$ where the space $S$ is either the torus $[-\pi, \pi]^d$, $\mathbb{Z}^d$,
or $\mathbb{R}^d$. Let $M$ denote the matrix with columns $\alpha _j,
j=1,...,k$, and let $r(A), r^*(A),$ denote the rank and dual rank of any set
$A$ of columns of $M$.

Let $f_j,\ j=1,\ldots,k$ be functions $f_j \in L_{p_j} (\mu(dx)),\ 1\leq
p_j\leq \infty$, defined on $S$, where $\mu(dx)$ is respectively normalized
Lebesgue measure, counting measure and Lebesgue measure.%

Let $z_j={\frac{1}{{p_j}}},\ j=1,\ldots,k, $ and $z=(z_1, ..., z_k)$. The
H\"older-Young-Brascamp-Lieb inequality ({GH}) will hold (with $K=K(z)<
\infty$) throughout the ``power counting polytopes" PCP defined respectively
by:

\medskip \noindent (c1) $\quad \displaystyle{\sum_{j\in A} z_j\leq r(A),\
\forall A}$

\noindent (c2) $\quad \displaystyle{\sum_{j\in A} z_j\geq r(M) -r(A^c),
\quad \Leftrightarrow \quad \sum_{j\in A} (1-z_j)\leq r^*(A) \ \forall A}$

\noindent (c3) $\quad \displaystyle{\sum_{j=1}^k z_j=m}$, and one of the
conditions (c1) or (c2) is satisfied.

Alternatively, the conditions (c1-c3) in the theorem are respectively
equivalent to:

\begin{enumerate}
\item ${\ z}=(z_1,\ldots,z_k)$ lies in the convex hull of the indicators of
the sets of independent columns of $M$, including the void set.

\item ${\ z}=(z_1,\ldots,z_k)$ lies in the convex hull of the indicators of
the sets of columns of $M$ which span its range.

\item ${\ z}=(z_1,\ldots,z_k)$ lies in the convex hull of the indicators of
the sets of columns of $M$ which form a basis.
\end{enumerate}

If, moreover, the matrix $M$ is unimodal, then the
H\"older-Young-Brascamp-Lieb inequality ({GH}) holds with constant $K=1$.
\end{thrm}

\textbf{Notes:} 1) Polytopes defined by the type of rank constraints
appearing in cases (c1-c2) are called ``polymatroids" (associated to $M$ and
$M^*$) -- see Welsh \cite{W}, 18.3, Theorem 1. The third polytope is the
intersection of the first two.

2) The first two cases of Theorem \ref{t:HYB} were obtained for unimodal
matrices in \cite{AB} and \cite{AT}, respectively.

3) Some further important issues, like the precise formula for $K(z)$, and
the nonhomogeneous extension where $l_j$ may be linear operators with
possibly different images, were resolved only recently -- see Lieb \cite{L},
Bennett, Carbery, Christ, Tao \cite{BCCT1}.

\textbf{Proof sketch:} By Edmonds theorem (see Welsh \cite{W}, 18.4, Theorem
1) the extremal points of the above polymatroids have only $0$ and $1$
coordinates for any matrix $M$. This fact leads to an easy proof, since at
the extremal points the result is immediate. For example, in the first and
third cases, the extremal points are in one to one correspondence with the
indicators of independent sets and bases $A$, respectively, and the constant
at such a point, by a change of variables, is seen to be:
\begin{equation*}
K_A=\sqrt{\det(A A^t)}
\end{equation*}

Since $K(z)$ is finite at the extremal points, and Riesz-Thorin
interpolation ensures the convexity of $\log(K(z))$, it follows that $K(z)$
will remain finite throughout the polytope generated by the indicators $1_A$.

\section{Wick products and Appell polynomials\label{s:appell}}

Let $W$ be a finite set and $Y_i$, $i\in W$ be a system of random variables.
Let
\begin{equation*}
Y^W=\prod_{i\in W}Y_i
\end{equation*}
denote the ordinary product, with $Y^\emptyset=1$, let $m^W = E \prod_{i\in
W}Y_i$ be the (mixed) moment, and let
\begin{equation*}
\chi(Y^W)=\chi(Y_i, i\in W)
\end{equation*}
denote the (mixed) cumulant of the variables $Y_i, i\in W$, defined
recursively as the solutions of the equations:
\begin{eqnarray}  \label{mts}
m^W = \sum_{\{V\} |- W}\chi(Y^{V_1})\cdots\chi(Y^{V_r}),
\end{eqnarray}
where the sum $\sum_{\{V\}|- W}$ is over all partitions $\{V\}=(V_1,%
\dots,V_r)$, $r\ge 1$ of the set $W$, and where $\chi(Y^\emptyset )=1$.

\textbf{Notes:} 1) The equation (\ref{mts}) is the formal power series
expression of the ``exponential relation $m=e^{\chi}$" between moments and
cumulants, viewed as functions on the lattice of subsets \cite{Niv}.

2) The inverse of the equation (\ref{mts}), the formal power series
expression of the ``logarithmic relation $\chi=\log(m)$" may also be
computed by:
\begin{equation*}
\chi(Y_1,\dots,Y_n)= {\frac{\partial ^T}{\partial z_1\dots \partial z_n}}
\log E\exp(\sum_{i=1}^Tz_jY_j)\Big|_ {z_1=\dots= z_n=0},
\end{equation*}
where the differentiation is interpreted formally if the moment generating
function does not exist.

\begin{dfntn}
The \textit{Wick products} $:Y^W:$ are defined as the solutions of the
recursion:
\begin{equation*}
Y^W = \sum_{U\subset W}:Y^U:E(Y^{W\setminus U})= \sum_{U\subset
W}:Y^U:\sum_{\{V\} |- W\setminus U}\chi(Y^{V_1})\cdots\chi(Y^{V_r}),
\end{equation*}
where the sum $\sum_{U\subset W}$ is taken over all subsets $U\subset W,$
including $U=\emptyset$, the sum $\sum_{\{V\}|- W\setminus U}$ is over all
partitions $\{V\}=(V_1,\dots,V_r)$, $r\ge 1$ of the set $W\backslash U$, and
the starting value is $:Y^\emptyset :=1$.
\end{dfntn}

\textbf{Notes:} 1) Inverting the recursion yields (\cite{surgailis:1983},
Proposition 1):
\begin{equation*}
:Y^W:= \sum_{U\subset W}Y^U\sum_{\{V\}|- W\setminus U}(-1)^r\chi(Y^{V_1})
\cdots\chi(Y^{V_r}),
\end{equation*}
as may be formally seen by replacing $m^{-1}$ by $e^{-\chi}$.

2) When some variables appear repeatedly, it is convenient to use the
notation
\begin{equation*}
:\underbrace{Y_{t_1},\ldots, Y_{t_1}}_{n_1},..., \underbrace {%
Y_{t_k},\ldots, Y_{t_k}}_{n_k}: = P_{n_1,...,n_k} (Y_{t_1},..., Y_{t_k})
\end{equation*}
(the indices in $P$ correspond to the number of times that the variables in
``$:\quad:$'' are repeated). The resulting multivariate polynomials $%
P_{n_1,...,n_k}$ are known as Appell polynomials. These polynomials are a
generalization of the Hermite polynomials, which are obtained if $Y_t$ are
Gaussian, and like them they play an important role in the limit theory of
quadratic forms of dependent variables (cf. \cite{surgailis:1983}, \cite%
{giraitis:surgailis:1986}, \cite{avram:taqqu:1987}).

4) The Appell polynomials may also be directly defined by "power-type"
recursions like:
\begin{eqnarray*}
&&\frac{\partial}{\partial \, x_j} P_{n_1,\ldots,n_k} (x_1,\ldots, x_k)=n_j
P_{n_1,\ldots,n_j -1, \ldots,n_k} (x_1,\ldots, x_k), \quad E
P_{n_1,\ldots,n_k} (X_1,\ldots, X_k)=0 \\
&&\forall n_j\geq 0, j=1,...,k, \ \sum_{j} n_j \geq 1, \\
&&P_{0,\ldots,0} (x_1,\ldots, x_k)=1.
\end{eqnarray*}

For example, when $m=n=1$, $P_{1,1}(X_{t},X_{s})=X_{t}X_{s}-{\mathbb{E}}%
X_{t}\,X_{s}$, and the bilinear form $Q_T(P_{1,1})$ is a weighted
periodogram with its expectation removed.

Note that the multivariate Appell polynomials can be defined by using
characteristic functions as well (see, e.g., \cite{surgailis:2000}).

\section{\label{s:cum}The diagram formula and the moments/cumulants of
sums/bilinear forms of Wick products}

\subsection{The cumulants diagram representation}

An important property of the Wick products is the existence of simple
combinatorial rules for calculation of the (mixed) cumulants, analogous to
the familiar diagrammatic formalism for the mixed cumulants of the Hermite
polynomials with respect to a Gaussian measure \cite{malyshev:1980}. Let us
assume that $W$ is a union of (disjoint) subsets $W_1,\dots, W_k$. If $%
(i,1),(i,2),\ldots,(i,n_i)$ represent the elements of the subset $W_i$, $%
i=1,\ldots,k$, then we can represent $W$ as a table consisting of rows $%
W_1,\dots, W_k$, as follows:
\begin{eqnarray}
\left (
\begin{array}{c}
(1,1),\ldots, (1,n_1) \\
\ldots\ldots\ldots \\
(k,1), \ldots, (k,n_k)%
\end{array}%
\right )=W.  \label{e:table}
\end{eqnarray}
By a \textit{diagram} $\gamma$ we mean a partition $\gamma=(V_1,\dots,V_r)$,
$r=1,2, \dots$ of the table $W$ into nonempty sets $V_i$ (the ``edges'' of
the diagram) such that $|V_i|\ge1$. We shall call the edge $V_i$ of the
diagram $\gamma$ \textit{flat}, if it is contained in one row of the table $%
W $; and \textit{free}, if it consists of one element, i.e. $|V_i|=1$. We
shall call the diagram \textit{connected}, if it does not split
the rows of the table $W $ into two or more disjoint subsets. We shall call
the diagram $\gamma=(V_1,\dots,V_r)$ \textit{Gaussian}, if $|V_1|=\dots=
|V_r|=2$. Suppose given a system of random variables $Y_{i,j}$ indexed by $%
(i,j)\in W$. Set for $V\subset W$,
\begin{equation*}
Y^V= \prod_{(i,j)\in V}Y_{i,j},\ \ \mbox{\rm and}\ \ \ :Y^V: \, = \,
:(Y_{i,j}, (i,j)\in V):\ .
\end{equation*}
For each diagram $\gamma =(V_1, \dots,V_r)$ we define the number
\begin{equation}
I_\gamma=\prod_{j=1}^r\chi (Y^{V_j}).  \label{e:Idef}
\end{equation}

\begin{prpstn}
\textrm{(cf. \cite{giraitis:surgailis:1986}, \cite{surgailis:1983})} \label%
{p:cum} Each of the numbers
\begin{eqnarray*}
& &(i)\quad EY^W= E(Y^{W_1}\dots Y^{W_k}), \\
& &(ii)\quad E(:Y^{W_1}:\dots :Y^{W_k}:), \\
& &(iii)\quad \chi(Y^{W_1},\dots, Y^{W_k}), \\
& &(iv)\quad \chi(:Y^{W_1}:,\dots , :Y^{W_k}:)
\end{eqnarray*}
\noindent is equal to
\begin{equation*}
\sum I_\gamma ,
\end{equation*}
where the sum is taken, respectively, over

(i) all diagrams,

(ii) all diagrams without flat edges,

(iii) all connected diagrams,

(iv) all connected diagrams without flat edges.

\noindent If $EY_{i,j}=0$ for all $(i,j)\in W$, then the diagrams in
(i)-(iv) have no singletons.
\end{prpstn}

\textbf{Notes:} 1) Part (i) is just the exponential relation between moments
and cumulants.

2) From part (ii) follows, for example, that $E:Y^W:=0$ (take $W=W_1$, then $%
W$ has only 1 row and all diagrams have flat edges).

\subsection{Multilinearity}

An important property of Wick products and of cumulants is their
multilinearity. For sums and bilinear forms
\begin{equation*}
S_{T}=S_{T}^m=\int_{I_T} P_m(X_{t})\nu(dt), \ \ \
Q_{T}=Q_{T}^{m,n}=\int_{I_T}\int_{I_T} \hat{b}(t-s)P_{m,n}(X_{t},
X_{s})\nu(dt)\nu(ds)
\end{equation*}
this implies that:
\begin{enumerate}
\item
\begin{eqnarray*}
&&\chi_k(S_T, ...,S_T)= \\
&& \qquad \int_{t_i\in I_T} \, \chi(:X_{t_{1,1}},\dots, X_{t_{1,m}}:,
\dots,:X_{t_{k,1}},\dots, X_{t_{k,m}}:) \prod_{i=1}^k \nu(dt_i),
\end{eqnarray*}
where the cumulant in the integral is taken for a table $W$ of $k$ rows $%
R_1...., R_k$, each containing the Wick product of $l$ variables identically
equal to $X_{t_k}$.
\item
\begin{eqnarray*}
&&\chi_k(Q_T, ...,Q_T)= \\
&& \qquad \int_{t_i, s_i \in I_T} \, \chi(:X_{t_{1,1}},\dots, X_{t_{1,m}},
\; X_{s_{1,1}},\dots, X_{s_{1,n}}:, \\
&& \qquad \qquad \qquad \dots,:X_{t_{k,1}},\dots, X_{t_{k,m}}, \;
X_{s_{k,1}},\dots, X_{s_{k,n}}:) \prod_{i=1}^k \hat{b}{(t_i-s_i)} \nu(dt_i)
\nu(d s_i),
\end{eqnarray*}
where the cumulant in the integral needs to be taken for a table $W$ of $k$
rows $R_1...., R_k$, each containing the Wick product of $m$ variables
identically equal to $X_{t_k}$ and of $n$ variables identically equal to $%
X_{s_k}$.
\end{enumerate}

\subsection{The cumulants of sums and quadratic forms of moving average
tables}

By part (iv) of Proposition \ref{p:cum}, applied to a table $W$ of $k$ rows $%
R_1...., R_k$, with $K=n_1 + ...n_k$ variables, and by the definition (\ref%
{e:Idef}) and of $I\gamma$, we find the following formula for the cumulants
of the Wick products of linear variables (\ref{nn7}):
\begin{equation}
\chi (:X_{t_{1,1}},\dots, X_{t_{1,n_1}}:,\dots,:X_{t_{k,1}},\dots,
X_{t_{k,n_k}}:) =\sum_{\gamma\in \Gamma (n_1,\ldots,n_k)}\kappa_{\gamma}
J_\gamma(\vec{t}) ,  \label{e:cum1}
\end{equation}
where $\Gamma (n_1,\ldots,n_k)$ denotes the set of all connected diagrams $%
\gamma=(V_1,\dots, V_r)$ without flat edges of the table $W$, $%
\kappa_{\gamma}=\chi_{|V_1|}(\xi_{I_1})\dots \chi_{|V_r|}(\xi_{I_1})$ and
\begin{eqnarray}  \label{e:cum2}
J_\gamma({t_1, ..., t_K}) &=& \prod_{j=1}^r J_{V_j}(t_{V_j}) \\
&=& \int_{s_1,\dots, s_r\in {\ I}} \prod_{j=1}^k \Big[\hat{a}%
(t_{j,1}-s_{j,1}) \hat{a}(t_{j,n_1}-s_{j,n_1})\dots  \notag \\
&& \quad \dots \hat{a}(t_{k,1}-s_{k,1})\dots \hat{a}(t_{k,n_k}-s_{k,n_k})%
\Big] \nu(ds_1),\dots, \nu(ds_r)  \notag \\
&=& \int_{\lambda_1, ...,\lambda_{K}} e^{i \sum_{j=1}^{K} t_j \; \lambda_j}
\; \prod_{i=1}^K a(\lambda_i) \; \prod_{j=1}^r \delta( \sum_{i \in V_j}
\lambda_i) \; \prod_{i=1}^{K} \mu(d \lambda_i)  \notag \\
&=& \int_{\lambda_1, ...,\lambda_{K}} e^{i \sum_{j=1}^{K} t_j \; \lambda_j}
\; \prod_{j=1}^r \left(f_{|V_j|}(\lambda_{j,1}, ..., \lambda_{j,|V_j|-1})
\delta( \sum_{i \in V_j} \lambda_i)\right) \; \prod_{i=1}^{K} \mu(d
\lambda_i),  \notag
\end{eqnarray}
where $s_{k,i}\equiv s_j$ if $(k,i) \in V_j$, $j=1,\dots,r$ and $%
\lambda_{j,i}\equiv \lambda_{i + \sum_{l=1}^{j-1}|V_l|}$ if $(j,i) \in V_j$,
$j=1,\dots,r$.

\subsection{The cumulants of sums and quadratic forms of moving average
processes.}

We will apply now the formula (\ref{e:cum2}) to compute the cumulants of $%
S_T^{(m)}, Q_T^{(m,n)}$. In this case, each row $j$ contains just one,
respectively two random variables.

It is easy to check that the variance of $S_T^{(2)}$ is:
\begin{equation*}
\chi_{2}(S_T^{(2)})= 2 \; \int_{\lambda_1,\lambda_2 \in S} \; {f}{(\lambda_1)%
} {f}{(\lambda_2)} \Delta_T( \lambda_1-\lambda_2) \Delta_T(
\lambda_2-\lambda_1) \prod_{e=1}^2 \mu({d \lambda_e}).
\end{equation*}
Note that there are only two possible diagrams on a table with two rows of
size $2$, and that they yield both a graph on two vertices (corresponding to
the rows), connected one to the other via two edges.

For another example, the third cumulant $\chi_{3}(S_T^{(2)})$ is a sum of
terms similar to:
\begin{equation*}
2^2 \int_{\lambda_1,\lambda_2,\lambda_3 \in S} \; {f}{(\lambda_1)} {f}{%
(\lambda_2)} {f}{(\lambda_3)} \Delta_T( \lambda_1-\lambda_2) \Delta_T(
\lambda_2-\lambda_3) \Delta_T( \lambda_3-\lambda_1) \prod_{e=1}^3 \mu({d
\lambda_e}).
\end{equation*}
This term comes from the $2^2$ diagrams in which the row $1$ is connected to
row $2$, $2$ to $3$ and $3$ to $1$.


For quadratic forms, a further application of part (iv) of Proposition \ref%
{p:cum} will decompose this as a sum of the form
\begin{equation*}
\sum_{\gamma\in \Gamma (n_1,\ldots,n_k)} \int_{t_i, s_i \in I_T}
R_\gamma(t_i, s_i) \prod_{i=1}^k \hat{b}_{t_i-s_i} d t_i d s_i,
\end{equation*}
where $\Gamma (n_1,\ldots,n_k)$ denotes the set of all connected diagrams $%
\gamma=(V_1,\dots, V_r)$ without flat edges of the table $W$ and $%
R_\gamma(t_i, s_i)$ denotes the product of the cumulants corresponding to
the partition sets of $\gamma$. This easy to check formula is also an
illustration of the diagram formula.

When $m=n=1$ and $k=2$, besides the Gaussian diagrams we have also one
diagram including all the four terms, which makes intervene the fourth order
cumulant of $X_t$.

When $m=n=1$, the Gaussian diagrams are all products of correlations and the
symmetry of $\hat{b}$ implies that all these $2^{k-1} (k-1)!$ terms are
equal. We get thus the well-known formula for the cumulants of discrete
Gaussian bilinear forms.

In general, we find decompositions as sums of certain ``Fej\'er graph
integrals", associated to specific graph structures.

The general structure of the intervening graphs for the cumulants of sums $%
S_T$ and quadratic forms $Q_T$ have been discussed in Section \ref{s:GS1}
(see Example \ref{exagraph}).

The following proposition is easy to check.

\begin{prpstn}
\label{p:cum-gauss} Let $X_{t},t\in I_{T}$ denote a stationary 
linear process given by (\ref{nn7}) with $d=1$. Then, the cumulants of the
sums and quadratic forms defined in (\ref{e:cum2}) are given respectively
by:
\begin{equation*}
\chi _{k,l}=\chi _{k}(S_{T}^{(m)},...,S_{T}^{(m)})=\sum_{\gamma \in \Gamma
(m,k)}\kappa _{\gamma }\;\sigma _{\gamma }(T)
\end{equation*}%
and
\begin{equation*}
\chi _{k,m,n}=\chi _{k}(Q_{T}^{(m,n)},...,Q_{T}^{(m,n)})=\sum_{\gamma \in
\Gamma (m,n,k)}\kappa _{\gamma }\;\tau _{\gamma }(T),
\end{equation*}%
where $\Delta _{T}(x)
$ is the Fej\'{e}r kernel, $\Gamma (l,k)$, $\Gamma (m,n,k)$ were defined
above, and
\begin{eqnarray}
\sigma _{\gamma }(T) &=&\int_{\vec{t}\in I_{T}^{k}}J_{\gamma }(\vec{t})dt
\label{cum1} \\
&=&\int_{\lambda _{1},...,\lambda _{K}}\prod_{j=1}^{k}\Delta
_{T}(\sum_{i=m(j-1)+1}^{m\,{j}}\;\lambda _{i})\;\prod_{i=1}^{K}a(\lambda
_{i})\;\prod_{j=1}^{r}\delta (\sum_{i\in V_{j}}\lambda
_{i})\;\prod_{i=1}^{K}d\lambda _{i},  \notag
\end{eqnarray}%
\begin{eqnarray}
\tau _{\gamma }(T) &=&\int_{\mu _{1},...,\mu _{k},\lambda _{1},...,\lambda
_{K},\lambda _{1}^{\prime },...,\lambda _{K^{\prime }}^{\prime }} \\
&&\qquad \prod_{j=1}^{k}\left[ \Delta _{T}\Big(\mu
_{j}+\sum_{i=m(j-1)+1}^{m\,{j}}\;\lambda _{i}\Big)\;\Delta _{T}\Big (-\mu
_{j}+\sum_{i=n(j-1)+1}^{n\,{j}}\;\lambda _{i}^{\prime }\Big)\;b(\mu _{j})%
\right]  \notag \\
&&\quad \times \prod_{i=1}^{K}a(\lambda _{i})\;\prod_{i=1}^{K^{\prime
}}a(\lambda _{i}^{\prime })\;\prod_{j=1}^{r}\delta (\sum_{i\in V_{j}}\lambda
_{i}+\sum_{i\in V_{j}}\lambda _{i}^{\prime })\;\;\prod_{i=1}^{K}d\lambda
_{i}\prod_{i=1}^{K^{\prime }}d\lambda _{i}^{\prime }\prod_{i=1}^{k}d\mu _{i}.
\notag  \label{cum2}
\end{eqnarray}
\end{prpstn}

These graph structures are simple enough to allow a quick evaluation of the
orders of magnitude $\alpha _M(z)$, via the corresponding graph-breaking
problems; for the case of bilinear forms we refer to Lemma 1 in \cite{a}.

For the case of sums, the domain of applicability of the CLT is $1-z_1 \geq
1/m$. We check now that at the extremal point $1-z_1 = 1/m$ we have
\begin{eqnarray*}  \label{e:profit}
\alpha _G(z_1) &=& \max_A p(A) \\
&=&\max_A [co(G-A) -\sum_{e \in A} (1- z_e)] \\
&=& \max_A [co(G-A) -|A| (1- z_1)] \\
&\leq& k/2,\ \forall G \in {\mathcal{G}}_k,
\end{eqnarray*}
where we interpret $p(A)$ as a ``profit," equal to the ``gain" $co(G-A)$
minus the ``cost" $\sum_{e \in A} (1- z_e)$. We thus need to show that at
the extremal point $1-z_1 = 1/m$,
\begin{equation*}
co(G-A) ]\leq |A|/m + k/2,\ \forall G \in {\mathcal{G}}_k.
\end{equation*}
Indeed, this inequality holds with equality for the ``total breaking" $A = {%
\mathcal{E}}$ (which contains $(k m)/2$ edges). It is also clear that no
other set of edges $A$ can achieve a bigger ``profit" $p(A)$ (defined in (%
\ref{e:profit})) than the total breaking, since for any other set $A$ which
leaves some vertex still attached to the others, the vertex could be
detached from the others with an increase of the number of components by $1$
and a cost no more than $m \frac{1}{m}$; thus the profit is nondecreasing
with respect to the number of vertices left unattached and thus the total
breaking achieves the maximum of $p(A)$.



\end{document}